\documentclass[11pt]{article}
 \usepackage{amssymb}
 \usepackage{epsfig}
 \usepackage{graphicx}
 \usepackage{pict2e}
 \newtheorem{Lemma}{Lemma}
 \newtheorem{Proposition}[Lemma]{Proposition}

 \newtheorem{OP}[Lemma]{Open Problem}

 \newtheorem{Corollary}[Lemma]{Corollary}

 \newcommand{\sfrac}[2]{{\textstyle\frac{#1}{#2}}}
 \newcommand{\eps}{\varepsilon}
 \newcommand{\Ints}{{\mathbb Z}}
 \newcommand{\Reals}{{\mathbb{R}}}
 
 \renewcommand{\AA}{\mbox{${\mathcal A}$}}
  
 \newcommand{\GG}{\mbox{${\mathbb{G}}$}}
 \newcommand{\RR}{\mbox{${\mathcal{R}}$}}

  \newcommand{\RRti}{\RR_{\mbox{{\footnotesize t-i} }}}
   \newcommand{\RRri}{\RR_{\mbox{{\footnotesize r-i} }}}
\newcommand{\TT}{\mbox{${\mathcal{T}}$}}
 
 \newcommand{\EE}{\mbox{${\mathcal{E}}$}}
 
      \newcommand{\rr}{\mbox{${\mathsf{r}}$}}
   \renewcommand{\ss}{\mbox{${\mathsf{s}}$}}
    \renewcommand{\SS}{\mbox{${\mathcal{S}}$}}
 
  \newcommand{\ed}{\ \stackrel{d}{=} \ }
   
  \newcommand{\bx}{{\mathbf x}}
   \newcommand{\bz}{{\mathbf z}}

  \newcommand{\proof}{\noindent {\bf Proof. }}
   \newcommand{\qed}{\ \ \rule{1ex}{1ex}}
    
   \newcommand{\area}{\mathrm{area}}
  \newcommand{\disc}{\mathrm{disc}}
   \newcommand{\circl}{\mathrm{circle}}
     \newcommand{\intensity}{\mathrm{intensity}}
   \newcommand{\len}{\mathrm{len}}
    \newcommand{\Leb}{\mathrm{Leb}}
   
     \newcommand{\cost}{\mathrm{cost}}
 \newcommand{\costbf}{\mbox{{\bf cost}}}
 \renewcommand{\time}{\mathrm{t}}
   \newcommand{\cst}{c_{\mbox{{\tiny ST}}}}
  
   \newcommand{\mcrit}{m_{\mbox{{\tiny crit}}}}

    \newcommand{\height}{\mathrm{height}}
  \newcommand{\peak}{\mathrm{peak}}
  \newcommand{\origin}{\mathbf{0}}
  \newcommand{\one}{\mathbf{1}}
   \renewcommand{\span}{\mathbf{span}}
  \newcommand{\Ex}{{\mathbb E}}
\renewcommand{\Pr}{\mathbb{P}}
 \begin{document}
  \title{Scale-invariant random spatial networks}
\author{David J. Aldous}

 \maketitle

\begin{abstract} 
Real-world road networks have an approximate scale-invariance property; 
can one devise mathematical models of random networks whose distributions are {\em exactly} invariant under Euclidean scaling?  
This requires working in the continuum plane. 
We introduce an axiomatization of a class of processes we call {\em scale-invariant random spatial networks}, whose 
 primitives are routes between 
each pair of points in the plane.  
We prove that one concrete model, based on minimum-time routes in a binary hierarchy of roads with different speed limits, satisfies the axioms, and note informally that two other constructions (based on Poisson line processes and on dynamic proximity graphs) are expected also to satisfy the axioms.  We initiate study of structure theory and summary statistics for general processes in this class.
\end{abstract}



\section{Introduction}
\label{sec-int}
Familiar web sites such as {\em Google maps} provide road maps on adjustable scale (zoom in or out) and a
suggested route between any two specified  addresses.  
Given $k$ addresses in a country, one could find the route for each of the 
${k \choose 2}$ pairs, and call the union of these routes the {\em subnetwork} 
(of the country's entire road network) spanning the $k$ points.

We abstract this idea by considering, for each pair of points $(z,z^\prime)$ in the plane,
a random route $\RR(z,z^\prime) = \RR(z^\prime,z)$ between $z$ and $z^\prime$.  
The collection of all routes (as $z$ and $z^\prime$ vary) defines what one might call a continuum random spatial network, an 
idea we  explain informally in this introduction 
(precise definitions will be given in section \ref{sec-TS}). 

In particular, for each finite set $(z_1,\ldots,z_k)$ of points we get a 
random network $\span(z_1,\ldots,z_k)$, the {\em spanning subnetwork}  linking the points, consisting 
of the union of the routes $\RR(z_i,z_j)$.
Mathematically natural structural properties we will impose on 
the distribution of such a process are 

\noindent
(i) translation and rotation invariance\\
(ii) scale-invariance. 
 
\vspace{0.03in}
\noindent 
For $0<c<\infty$ the scaling map $\sigma_c: \Reals^2 \to \Reals^2$  takes $z$ to
$cz$; we emphasize that (ii) means ``naive Euclidean scaling", i.e. invariance under the action of $\sigma_c$, not any notion of ``scaling exponent".  
For instance, scale-invariance implies that the route-length $D_r$ between points 
at (Euclidean) distance $r$ apart must scale as $D_r \ed r D_1$, where of course 
$ 1 \le D_1 \le \infty $.
The setup so far does not exclude the possibility that routes are 
fractal, with infinite length, and such cases do in fact arise naturally in the 
tree-like models of section \ref{sec-tree}.  
But, envisaging road networks rather than some other physical structure, we 
restrict attention to the case 
$\Ex D_1 < \infty$.  
There is a rather trivial example, the {\em complete network} in which each 
$\RR(z_1,z_2)$ is the straight line segment from $z_1$ to $z_2$, 
but the assumption ``$\ell < \infty$" below will exclude this example.
 
Much of our study involves 
{\em sampled spanning subnetworks} $\SS(\lambda)$, as follows.
Write $\Xi(\lambda)$ for a Poisson point process in $\Reals^2$ of  intensity $\lambda$, 
independent of the network.
Then the points $\xi$ of $\Xi(\lambda)$, together with the routes $\RR(\xi,\xi^\prime)$ for
each pair of such points, 
form a random subnetwork we denote by $\SS(\lambda)$.  
The distribution of $\SS(\lambda)$ inherits the properties 
of translation- and rotation-invariance, and a form of scale-invariance described at (\ref{SS-def-4}).
In particular we can define a constant $0 < \ell \le \infty$ by 
\[ \ell = \mbox{mean length-per-unit-area of } \SS(1) \]
(where ``mean length-per-unit-area " is formalized by {\em edge-intensity} at (\ref{intensity-def})).
In section \ref{sec-Steiner} we note a crude lower bound $\ell \ge \frac{1}{4}$. 
We impose the property 
\[
 \ell  < \infty .
\]
Regard $\ell$ as ``normalized network length", for the purpose of comparing different networks.

Everything mentioned so far makes sense when only 
finite-dimensional distributions $\RR(z_i,z_j)$ are specified. 
A first context in which we want to consider a
process over the whole continuum concerns the following convenient abstraction of the notion of ``major road". 
Write $\RR_{(1)}(z_1,z_2)$ for the part of the route $\RR(z_1,z_2)$ that is at distance $\geq 1$ from each of $z_1$ and $z_2$.
Conceptually, we want to study an edge-process $\EE$ viewed as the union of $\RR_{(1)}(z_1,z_2)$ over all pairs $(z_1,z_2)$.  
To formalize this directly would require 
some notion of
``regularity" for a realization, for instance some notion of a.e. continuity of routes  $\RR(z_1,z_2)$
as $z_1$ and $z_2$ vary. 
But we can avoid this complication by first considering only $z_1, z_2$ in
$\Xi(\lambda)$ 
and then letting $\lambda \to \infty$.  
After defining $\EE$ in this way, we can define 
\[
p(1)  := \mbox{mean length-per-unit-area of } \EE 
\]
and impose the requirement
\[
p(1) < \infty . 
\]
If a process of random routes $\RR(z,z^\prime)$  satisfies the properties we have described 
(as stated precisely in section \ref{sec-TS}), 
then we will call it a {\em scale-invariant random spatial network} (SIRSN).  
As the choice of name suggests, it is the scale-invariance that makes such processes of mathematical interest; 
in section \ref{sec-visualizing} we briefly discuss its plausibility for real-world networks.

We do not know any closely related previous work.   
We will discuss one related area of theory  (discrete random spatial networks; section \ref{sec-spatial}) 
and one area of application 
(fast algorithms for shortest routes:; section \ref{sec-intro-SPA}).
Several more distantly related topics are mentioned in section \ref{sec -ORL}.

\subsection{Outline of paper}
The purpose of this paper is to initiate study of SIRSNs,  with three emphases.  
First, we give a careful formulation of an axiomatic setup for SIRSNs, 
with discussion of possible alternatives (section \ref{sec-2}).  
Second, it is not obvious that  SIRSNs exist at all!  
We give details of one construction in  section \ref{sec-constr1}.
That construction envisages a square lattice of freeways, with ``speed level $j$" freeways spaced $2^j$ apart, and 
the routes are the minimum time paths. 
Being based on the discrete lattice makes some estimates technically straightforward, but  completing the details of proof requires surprisingly intricate arguments.
This construction is somewhat artificial in not naturally having all the desired invariance properties, so these need to be forced by external randomization. 
We briefly mention two other constructions 
(in section \ref{sec-constr2} based on a weighted Poisson line process representing the different-level freeways, and in section \ref{sec-constr3} based on a dynamic construction of random points and roads added accoding to a deterministic rule)
which intuitively seem more natural but for which we have been unable to complete all the details of a proof.

Third, in sections \ref{sec-properties} - \ref{sec-stronger} we begin developing some general theory from the axiomatic
setup.
Of course {\em scale-invariance} is a rather weak assumption, 
loosely analogous to {\em stationarity} for a stochastic process, 
so one cannot expect sharp results holding throughout this general class of process.  
Our general results might be termed ``structure theory" and concern 
existence and uniqueness issues for singly- and doubly-infinite geodesics, 
continuity of routes $\RR(z_1,z_2)$ as a function of $(z_1,z_2)$, 
numbers of routes connecting disjoint subsets,
and bounds on the parameters $\Ex D_1, \ell, p(1)$.
One feature worth emphasis is that 
(very loosely analogous to {\em entropy rate} for a stationary process) 
the quantity $p(1)$ is a non-obvious statistic of a SIRSN, but turns out to be 
central in the foundational setup, in the structure theory, and in 
 conceptual interpretation as a model for road networks.    
The latter is best illustrated by the ``algorithms" story in sections \ref{sec-intro-SPA} and \ref{sec-transit}.

Being a new topic there are numerous open problems, both conceptual and technical, 
stated in a final discussion section \ref{sec-final-disc}.

Before starting technical material,  sections \ref{sec-spatial} - \ref{sec-visualizing}  give further verbal discussion of background 
to the topic.

\subsection{Discrete spatial networks}
\label{sec-spatial}
Traditional models of (deterministic or random) spatial networks start with a {\em discrete} set of points
and then assign linking edges via some rule, 
e.g.  the random geometric graph \cite{penrose-RGG} or proximity graphs \cite{jaromczyk}, surveyed in  \cite{me-spatial-1}.
One specific motivation for the present work was as a second attempt to resolve a paradox -- 
more accurately, an unwelcome feature of a naive model -- in the discrete setting, observed in \cite{me116}.
In studying the trade-off between total network length and the effectiveness of  
a network in providing short routes between discrete cities, one's first thought might be to measure the
latter by
the average, over {\em all} pairs $(x,y)$, of the ratio
\[
 \mbox{(route-length for $x$ to $y$)}/\mbox{(Euclidean distance from $x$ to $y$)} 
\]
instead of averaging over pairs at  Euclidean distance $\approx r$ to get our $\Ex D_r$.
But it turns out that (in the $n \to \infty$ limit of a network on $n$ points) one
can make this ratio tend to $1$ 
for a network whose length is only $1 +o(1)$ times the length of the Steiner tree,
by simply 
superimposing on the Steiner tree a sparse Poisson line process.  
Such ``theoretically optimal" networks are completely unrealistic, so there must be something wrong with the
optimization criteria.  
What's wrong is that the networks are ineffective for small $r$.
One way to get a non-trivial tradeoff in the $n \to \infty$ limit was described in 
\cite{me-spatial-1}: 
using the statistic $\max_r  r^{-1} \Ex D_r$ as the measure of effectiveness leads to more realistic-looking networks.
In the discrete setting a network model cannot be precisely scale-invariant, but 
such  considerations prompted investigation of continuum models which are assumed to be scale-invariant, so that 
$r^{-1} \Ex D_r$ is constant.

We emphasize that our networks involve roads at definite positions in the plane.  There is substantial recent literature, discussed in \cite{itai-planar}, involving 
quite different notions of random planar networks, based on identifying topologically equivalent 
networks.

\subsection{Visualizing spanning subnetworks}
Both  construction and analysis of general SIRSNs are based on
studying subnetworks $\span(z_1,\ldots,z_k)$ on fixed or (most often) random points.
It is helpful to visualize what  subnetworks look like -- see Figure 1.

\setlength{\unitlength}{0.06in}
\begin{picture}(40,40)(-9,0)

\put(5.9,35.9){\circle*{0.9}}
\put(5.9,35.9){\line(-1,-1){0.9}}
\put(2,34){\line(1,1){2}}
\put(5,35){\line(-2,1){1}}
\put(4,35.5){\line(0,1){0.5}}
\put(35,5){\circle*{0.9}}
\put(35,5){\line(-1,0){1}}
\put(34,5){\line(-1,-1){1}}
\put(33,4){\line(-3,1){3}}
\put(30,5){\line(-1,1){27.5}}
\put(2.5,32.5){\line(-1,3){0.5}}
\put(32.7,32.7){\circle*{0.9}}
\put(32,32){\line(1,1){0.7}}
\put(32,32){\line(-1,0){1}}
\put(31,32){\line(-1,2){1}}
\put(30,34){\line(-6,-1){18}}
\put(12,31){\line(-1,1){4.5}}
\put(7.5,35.5){\line(-2,-1){2.5}}
\put(5,34.25){\line(0,1){0.75}}
\put(31,32){\line(-1,-1){2}}
\put(29,30){\line(-1,-3){6}}
\put(14,13.2){\circle*{0.9}}
\put(14,13.2){\line(0,1){0.8}}
\put(14,14){\line(3,-1){3}}
\put(17,13){\line(3,2){3}}
\put(21,14){\line(1,0){2.67}}
\put(14,14){\line(-2,1){2}}
\put(12,15){\line(0,1){2}}
\put(12,17){\line(2,1){4}}
\put(32.3,20.3){\circle*{0.9}}
\put(33,21){\line(-1,-1){0.7}}
\put(33,21){\line(-1,1){5.2}}
\put(31,23){\line(-5,-2){3}}
\put(28,21.8){\line(-80,32){8}}
\put(20,25){\line(-1,5){1.42}}
\put(4.6,9.6){\circle*{0.9}}
\put(4,9){\line(1,1){0.6}}
\put(4,9){\line(-1,1){1}}
\put(3,10){\line(0,1){2.5}}
\put(3,12.5){\line(2,1){24.2}}
\put(4,9){\line(3,-1){3}}
\put(7,8){\line(4,1){8}}
\put(15,10){\line(3,-1){15}}
\put(8.3,15.06){\line(-1,3){5.8}}
\put(8.5,21.5){\circle*{0.9}}
\put(8,21){\line(1,1){0.5}}
\put(8,21){\line(-2,1){2}}
\put(8,21){\line(10,-18){1}}
\put(9,19.2){\line(5,-1){8.5}}
\put(13,17.5){\line(0,1){0.9}}
\end{picture}

{\bf Figure 1.}
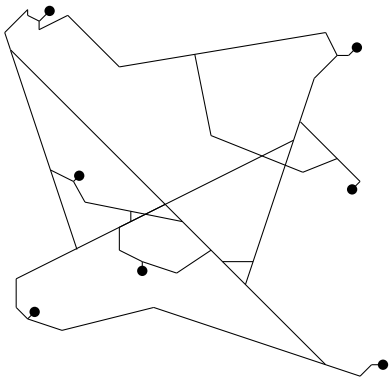
{\small Schematic for the subnetwork of a SIRSN on 7 points $\bullet$}

\medskip
\noindent
The qualitative appearance of Figure 1 is  quite different from that of  familiar spatial networks 
mentioned above, 
based on a discrete set of points, which could be viewed as 
 abstractions of an inter-city road network, with cities as points.  
In contrast, we are abstracting the idea of the points $\bullet$ being  individual street addresses a long way apart.
The real-world route between two such street addresses will typically consist, in the middle, of roughly straight freeway segments but,
nearing an endpoint, of a more jagged trajectory of shorter segments of lower-capacity roads; our setup and 
Proposition \ref{Pjagged} imply the same behavior in our model.

\subsection{Very fast shortest path algorithms}
\label{sec-intro-SPA}
There is an interesting connection with the ``shortest path algorithms" literature. 
Online mapping services and GPS devices require very quick computations of shortest routes. 
In this context,   
the U.S road network is represented as a graph on about 15 million street intersections  
(vertices) with edges (road segments) marked by distance 
(or typical driving time), and a given street address is recognized as being between two specific street intersections.
Given a pair of (starting and destination) points, one wants to compute the shortest route.
Neither of the two extremes -- pre-compute and store the routes for all possible pairs; 
or use a classical Dijkstra-style algorithm for a given pair without any preprocessing -- is practical.  
Bast et al (see \cite{bast-science} for an outline)  find a set of about 10,000 intersections 
(which they call {\em transit nodes}) 
with the property that, unless the start and destination points are close, the shortest route goes via 
some transit node near the start and 
some transit node near the destination. 
Given such a set, one can pre-compute shortest routes and route-lengths between each pair of transit nodes; 
then answer a query by using the classical algorithm to calculate the route lengths from starting 
(and from destination) point to each nearby transit node, and finally minimizing over pairs of such transit nodes.

This idea is actually used commercially (and patented). 
Mathematical discussion was initiated by 
Abraham et al \cite{goldberg10}, who 
introduced the notion of {\em highway dimension}, 
defined as the smallest integer $h$ such that 
for every $r$ and every ball of radius $4r$, 
there exists a set of $h$ vertices such that 
every shortest route of length $>r$ within the ball passes through some vertex in the set. 
They discuss several algorithms whose performance can be analyzed in terms of highway dimension,  
and devise a particular model (a dynamic spanner construction on vertices given by an adversary) 
designed to have bounded highway dimension.

Now saying one can find $h$ independent of $r$ is a form of approximate scale-invariance, so the empirical fact that
one can find transit nodes in the real-world road networks is a weak form of empirical scale-invariance.
Within our model where precise scale-invariance is assumed, we can derive quantitative estimates relating 
to transit nodes -- see section \ref{sec-transit}.

Incidently, the way we define edge-processes $\EE = \EE(\lambda,r)$ in terms of routes
(mentioned in the Introduction and defined in section \ref{sec-TS}) is closely related to
the notion of {\em reach} in the algorithmic literature \cite{Gutman04}.

\subsection{Visualizing scale-invariance}
\label{sec-visualizing}
Visualizing a photo of a road, scale-invariance seems implausible, because it implies existence of roads of 
arbitrarily large and arbitrarily small ``sizes", however one interprets ``size".   
But scale-invariance is not referring to the physical roads but to the process of ``shortest routes", as  in 
the discussion above.  
Figure 2  illustrates one aspect of scale-invariance. 
There is some number of crossing places (over the line) used by routes from one square to the other square. 
In our model, scale-invariance implies that the mean number of such crossings 
does not depend on the scale of the map.  
One could test this as a prediction about real-world road networks.

\setlength{\unitlength}{0.6in}
\begin{picture}(4,2,4)(-1.2,-0.6)
\put(0,0){\line(1,0){1}}
\put(0,0){\line(0,1){1}}
\put(1,1){\line(-1,0){1}}
\put(1,1){\line(0,-1){1}}
\put(3,0){\line(1,0){1}}
\put(3,0){\line(0,1){1}}
\put(4,1){\line(-1,0){1}}
\put(4,1){\line(0,-1){1}}
\put(2.0,-0.5){\line(0,1){2}}
\put(1.6,0.8){\line(3,1){0.8}}
\put(1.6,-0.1){\line(5,1){0.8}}
\put(1.6,0.5){\line(6,-1){0.8}}
\end{picture}

{\bf Figure 2.} 
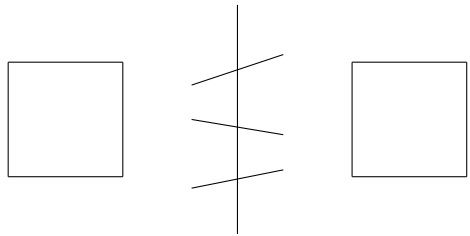
{\small Schematic for long-distance routes.}

\medskip
\noindent
As another empirical aspect of scale-invariance,
 \cite{PhysRevE.73.026130} studied proportions of route-length, within distance-$r$ routes,  spent on the $i$'th longest road segment in the route 
(identifying roads by their highway number designation) and observe that in the U.S. the averages 
of these ordered proportions are around 
$(0.40, 0.20, 0.13, 0.08, 0.05)$ as $r$ varies over a range of medium to large distances.  
Again, in our models (identifying roads as straight segments) scale-invariance implies there is some vector 
of expected proportions that is precisely independent of $r$.

\section{Technical setup}
\label{sec-2}

In formulating an axiomatic setup there are several alternative choices one could make. 
In section \ref{sec-TS} we state concisely the choices we made; 
section \ref{sec-TS-2}  discusses alternatives, reasons for choices, and immediate consequences or non-consequences of the setup.

\subsection{Stochastic geometry background}
We quote a  fundamental identity from stochastic geometry (see \cite{MR895588} Chapter 8).
Let $\EE$ be an {\em edge process} -- for our purposes, a union of line segments --  whose distribution is invariant under translation and rotation.  
Then $\EE$ has an {\em edge-intensity}, a constant $\iota = \intensity(\EE) \in [0,\infty]$ such that 
\begin{equation}
 E (\mbox{length of } \EE \cap A) = \iota \times \area(A), \quad A \subset \Reals^2 . 
\label{intensity-def}
\end{equation}
Moreover the positions and angles at which $\EE$ intersects the $x$-axis 
(and hence any other line) are such that 
\begin{equation}
\mbox{mean number intersections per unit length} = 2 \pi^{-1} \times \intensity(\EE) 
\label{intersection-identity}
\end{equation}
and the random angle $\Theta \in (0,\pi)$ of a typical intersection has density
\begin{equation}
f_\Theta(\theta) = \sfrac{1}{2} \sin \theta . 
\label{angle}
\end{equation}

\subsection{Definitions}
\label{sec-TS}
Here we organize the setup via four aspects.

\paragraph{Some notation.}  
$\origin$ denotes the origin; 
$\disc(z, r)$ and $\circl(z, r)$ denote the closed disc and the circle centered at $z$.

\paragraph{Aspect 1.  Allowed routes and route-compatability.}  
 Define a {\em  jagged route} between two points $z, z^\prime$ of $\Reals^2$ 
to  consist of straight line segments 
between successive points $(z_i, - \infty < i < \infty) $ with $\lim_{i \to - \infty} z_i = z$ and $\lim_{i \to \infty} z_i = z^\prime$, and such that the total length 
$\sum_{i = -\infty}^\infty |z_i - z_{i-1}|$ is finite.  
A {\em feasible route} is either a jagged route or the variant
with a finite or semi-infinite set of successive line segments; 
we further require that the route be non-self-intersecting.
Write  $\mathsf{r}(z,z^\prime)$ for a feasible route, which from now on we will just call {\em route}.  
We envisage a route $\rr(z,z^\prime)$ as a one-dimensional subset of $\Reals^2$, equipped with a label indicating it is the route from 
$z$ to $z^\prime$.   The route $\rr(z^\prime,z)$ is always the reversal of $\rr(z,z^\prime)$.

When we have a collection of  routes, we require the following 
{\em pairwise compatability} property. 

\vspace*{-0.1in}
\[ \mbox{  If two routes $\mathsf{r}(z_1,z_j), \ \mathsf{r}(z^\prime_1,z^\prime_2)$ meet at two points then the routes } \] 

\vspace*{-0.3in}
\begin{equation}
\mbox{coincide on the subroute between the two meeting points.} 
\end{equation}

\paragraph{Aspect 2.    Subnetworks on locally finite configurations.} 
Given a locally finite configuration of points $(z_i)$ in the plane, and routes $\rr(z_i,z_j)$ satisfying the pairwise compatability property, 
write $\ss$ for the union of all these routes.  
If $\ss$ has the ``finite length in bounded regions" property 
\begin{equation}
\len(\ss \cap \disc(\origin,r) ) < \infty \mbox{ for each } r < \infty 
\label{BLBA}
\end{equation}
then call $\ss$ a {\em feasible subnetwork}. 
Here ``$\len$' denotes ``length".
Formally $\ss$ consists of the vertex set 
$(z_i)$, an edge set which is the union of the edge sets comprising each 
$\rr(z_i,z_j)$, and marks on edges to indicate which routes they are in. 
Inclusion $\ss(1) \subseteq \ss(2)$ means that $\ss(2)$ can be obtained from 
$\ss(1)$ by adding extra vertices and associated routes.

As outlined in section \ref{sec-TS-2} there is a natural $\sigma$-field that makes 
the set of all feasible subnetworks into a measurable space, 
so it makes sense below to talk about random feasible subnetworks.

\paragraph{Aspect 3.  Desired distributional properties of subnetworks.}
The precise definition of the class of processes we shall study uses 
``finite-dimensional distributions" (FDDs), as follows.
Given a finite set $z_1,\ldots,z_k$ 
let $\mu_{z_1,\ldots,z_k}$ be the distribution of a random 
feasible subnetwork $\span(z_1,\ldots,z_k)$ on $z_1,\ldots,z_k$.
Suppose a family (indexed by all finite sets) of FDDs satisfies 
\begin{eqnarray}
&&\mbox{the natural consistency condition} \label{fdd-1} \\
&&\mbox{invariance under translation and rotation} \label{fdd-2} \\
&&\mbox{invariance under scaling.} \label{fdd-3} 
\end{eqnarray}
To be precise about (\ref{fdd-3}), recall that the scaling map $\sigma_c: \Reals^2 \to \Reals^2$  takes $z$ to $cz$.  
Then the action of $\sigma_c$ on $\span(z_1,\ldots,z_k)$ gives a random subnetwork
whose distribution 
equals the distribution of $\span(\sigma_c z_1,\ldots,\sigma_c z_k)$. 

Appealing to the Kolmogorov extension theorem, we can associate 
with such a family a process of routes $\RR(z_1,z_2)$, for each pair $z_1, z_2$ in
$\Reals^2$, 
though for a process defined in that way we can only discuss properties determined by FDDs.

As mentioned earlier, much of our study involves 
{\em sampled spanning subnetworks}, as follows.
For each $0<\lambda<\infty$ let 
 $\Xi(\lambda)$ be a Poisson point process of intensity
$\lambda$   (we sometimes call this {\em 
point-intensity} to distinguish from edge-intensity at (\ref{intensity-def})). 
Make a process  $(\Xi(\lambda), \ 0<\lambda<\infty)$ by
coupling in the natural way 
(take a space-time Poisson point process and let $\Xi(\lambda)$ be the positions of
points arriving during time $[0,\lambda]$). 
Taking Poisson points independent of the  process of routes,  we can  define $\SS(\lambda)$ as the subnetwork of routes $\RR(\xi, \xi^\prime) $ for 
pairs $\xi, \xi^\prime$ in $\Xi(\lambda)$.
We want the resulting processes $\SS(\lambda)$ to have  the following properties. 

{\small 
\begin{eqnarray}
&&\mbox{for each $\lambda$, $\SS(\lambda)$ is a random feasible subnetwork on vertex-set $\Xi(\lambda)$} \label{SS-def-1}\\
&&\mbox{for each $\lambda$,  $\SS(\lambda)$ has translation- and rotation-invariant distribution }  \label{SS-def-2}\\
&&\mbox{$\SS(\lambda_1) \subseteq \SS(\lambda_2)$ for $\lambda_1 < \lambda_2$}\label{SS-def-3} \\
&&\mbox{applying $\sigma_c$ to $\SS(\lambda)$ gives a network distributed as $\SS(c^{-2} \lambda)$} \label{SS-def-4} .
\end{eqnarray}
}
For (\ref{SS-def-4}), recall that 
applying $\sigma_c$ to $\Xi(\lambda)$ gives a point process distributed as $\Xi(c^{-2} \lambda)$.

We omit full measure-theoretic details of the construction of $\SS(\lambda)$, and just 
point out what extra conditions are needed to obtain properties 
(\ref{SS-def-1} - \ref{SS-def-4}).
First, we need to impose the technical condition
\begin{equation}
\mbox{the map $(z_1,\ldots,z_k) \to \mu_{z_1,\ldots,z_k} $ is measurable} 
\label{fdd-4}
\end{equation} 
to ensure that $\SS(\lambda)$ is measurable.  
Second, part of the ``feasible" assertion in (\ref{SS-def-1}) is the  ``finite length in bounded regions" property (\ref{BLBA}),
and this property for $\SS(\lambda)$ cannot be a consequence of assumptions on FDDs only, so we need 
\begin{equation}
\len(\SS(\lambda) \cap \disc(\origin,r) ) < \infty \ \mbox{ a.s. for each } r < \infty 
\label{fdd-5} 
\end{equation}
and this will follow from the stronger assumption (\ref{ell-finite}) below.

\paragraph{Aspect 4.  Final definition of a SIRSN.}
To summarize the above:
given a process of routes $\RR(z_1,z_2)$ with FDDs 
satisfying (\ref{fdd-1} - \ref{fdd-3}, \ref{fdd-4}), 
we can define the process of sampled subnetworks $(\SS(\lambda), 0 < \lambda < \infty)$ 
which, if (\ref{fdd-5}) holds, will have properties (\ref{SS-def-1} - \ref{SS-def-4}).
Finally, we define a SIRSN as a process 
(denoted by the routes $\RR(z_1,z_2)$ or by the sampled subnetworks
$(\SS(\lambda), 0 < \lambda < \infty)$) satisfying these assumptions 
(\ref{fdd-1} - \ref{fdd-3}, \ref{fdd-4}) 
and also satisfying the extra conditions (\ref{D1},\ref{p1intro}) below.  
These extra conditions merely repeat and formalize  the requirements, stated in the 
introduction, that certain statistics be finite.  
As noted above, these assumptions imply that  (\ref{SS-def-1} - \ref{SS-def-4}) hold.

Write  $\one = (1,0)$ and $D_1: = \len \  \RR(\origin, \one)$. 
So $D_1$ represents route-length between points at distance $1$ apart. 
Our definition of {\em feasible} route implies 
$1 \le D_1 < \infty$ a.s., and we impose the requirement
\begin{equation}
 1 < \Ex D_1 < \infty .
\label{D1}
\end{equation} 
Next, our definition of {\em feasible subnetwork} implies that $\SS(1)$ must have 
a.s. finite length in a bounded region.  We impose the stronger requirement 
of finite {\em expected} length.  
In terms of the edge-intensity (\ref{intensity-def}), we require
\begin{equation}
 \ell := \intensity( \SS(1)) < \infty . 
\label{ell-finite}
\end{equation}
Finally, we define 
\begin{equation}
\EE(\lambda, r) := \bigcup_{\xi, \xi^\prime \in \Xi(\lambda)} \ \ 
\RR(\xi,\xi^\prime) \setminus (\disc(\xi,r) \cup \disc(\xi^\prime,r)) 
\label{EE-def}
\end{equation}
and edge-intensities
\begin{eqnarray}
p(\lambda, r) &:=& \intensity(\EE(\lambda,r)) \label{plar-def}\\ 
p(r) &:=& \lim_{\lambda \to \infty} p(\lambda,r) \label{pr-def}
\end{eqnarray}
and impose the requirement
\begin{equation}
p(1)   < \infty 
\label{p1intro}
\end{equation} 
whose significance is discussed in the next section.
Lemma \ref{Lellp} will show that (\ref{p1intro}) implies (\ref{ell-finite}).
If we do not require (\ref{p1intro})  but instead require (\ref{ell-finite}), 
call the process a {\em weak} SIRSN.

\subsection{Discussion of technical setup}
\label{sec-TS-2}

\paragraph{Aspect 1.  Allowed routes and route-compatability.}  
Because we want routes to have a well-defined lengths, a minimum assumption would be that routes 
are rectifiable curves.  
We have assumed  ``feasible routes"  in order to simplify notation.
We believe that the theory would be essentially unchanged if instead one allowed rectifiable curves, 
as in the (quite different) theory mentioned in section \ref{sec-MK}. 

It turns out (section \ref{sec-jagged}) that realizations of our models always have 
jagged routes. 
A consequence is that (as in Figure 1) a route 
$\RR(\xi,\xi^\prime)$ between two points of $\SS(\lambda)$ does not pass through any third point $\xi^{\prime \prime}$ of $\SS(\lambda)$.  
This prompts the precise definition of {\em geodesic} below.

The route-compatability property is a property that would hold if routes were defined as minimum-cost paths, for some reasonable notion of
``cost".  
Note that our formal setup does not require routes to be minimum-cost in any explicit sense.  

\paragraph{Aspect 2.    Subnetworks on locally finite configurations.} 
Here are some properties of a fixed feasible subnetwork.
\begin{Lemma}
\label{Lss1}
Let $\ss$ be a feasible subnetwork on a locally finite, infinite configuration
$(z_i)$.\\
(i) The set 
$ \{ \rr(z_i,z_j) \cap \disc(z, r)\}_{i,j}$ of sub-routes appearing as intersections 
of some route with a fixed disc $\disc(z, r)$ contains only finitely many distinct
(non-identical) 
sub-routes. \\
(ii) For each $i$ and each sequence $(z_j)$ with $|z_j| \to \infty$ there is a
subsequence 
$z^\prime_k = z_{j(k)}$ and a semi-infinite path $\pi$ from $z_i$ in $\ss$ such
that, for each $r > 0$, 
\[ \rr(z_i,z^\prime_k) \cap \disc(z_i, r) = \pi \cap \disc(z_i, r) 
\mbox{ for all large } k . \]
\end{Lemma}
\proof
(ii) follows from (i) by a compactness argument.  
To outline (i), if false then 
(by the route-compatability property) 
the subroutes must meet the disc boundary at an infinite number of distinct points, 
and then 
(again by the route-compatability property) 
their extensions 
must meet the  boundary of a slightly larger disc at an infinite number of distinct
points, 
implying infinite length and contradicting the ``finite length in bounded regions" 
property (\ref{BLBA}) of $\ss$.
\qed

Note that Lemma \ref{Lss1} is implicitly about compactness in a topology on the space of paths within a given subnetwork $\ss$.
This is quite different from the topology of the space of all subnetworks, mentioned later.

\paragraph{Terminology: paths, routes and geodesics.}
  A {\em path} in $\ss$ has its usual network meaning.  
Typically there will be many paths between $z_i$ and $z_j$, but 
(as part of the structure of a {\em feasible subnetwork}) 
one is distinguished as the {\em route} $\rr(z_i,z_j)$.  
So a route is a path; and a path may or may not be part of one or more routes.  
 A {\em singly infinite geodesic} in $\ss$ from $z_i$  is an infinite path, starting from $z_i$, 
 such that any finite portion of the path is a subroute of the route $\rr(z_i,z_k)$ for some $z_k$. 
 So Lemma \ref{Lss1}(ii) says that there always exists at least one singly infinite geodesic from $z_i$. 
 A typical point $\epsilon$ along a route $\rr(z_i,z_j)$ will sometimes be called a {\em path element} to distinguish it from the endpoints.

Now write $\mathfrak{S}$ for the set of all feasible subnetworks $\ss$ on all locally finite 
configurations $\bx = (x_j)$.  
It is natural to want to regard $\SS(\lambda)$ as a random element of $\mathfrak{S}$, 
which requires specifying a $\sigma$-field on $\mathfrak{S}$, 
and as traditional we can do this by specifying 
a complete separable metric space structure on $\mathfrak{S}$ 
and using the Borel $\sigma$-field.

We outline a ``natural" topology in an appendix.
In this paper the topology plays no explicit role, 
but one can imagine  developments where it does -- one can imagine constructions using weak convergence, 
for instance, and compactness issues would be key to a proof of the existence part 
 of Open Problem \ref{OP2}.  
However, it might be better to develop such theory within a framework where routes are allowed to be rectifiable curves.

\paragraph{Aspect 3.  Desired distributional properties of subnetworks.}

The scale-invariance property (\ref{SS-def-4})
\[ \mbox{applying $\sigma_c$ to $\SS(\lambda)$ gives a network distributed as $\SS(c^{-2} \lambda)$} \] 
is what gives SIRSNs a mathematically interesting structure, and almost all our general results 
 in sections \ref{sec-properties} and \ref{sec-stronger} rely on scale-invariance.  
 To indicate how it is used, 
define $\ell(\lambda)$ analogously to (\ref{ell-finite}): 
\begin{equation}
\ell(\lambda) :=   \intensity( \SS(\lambda)) 
\label{ell-lam-def}
\end{equation} 
so $\ell(1) = \ell$.  
Then there is a scaling relation
\begin{equation}
\ell (\lambda) = \lambda^{1/2} \ell, \quad 0 < \lambda < \infty . 
\label{ell-scale}
\end{equation}
To derive this relation, consider
the scaling map $\sigma_{\lambda^{-1/2}}$ that takes 
$\SS(1)$ to $\SS(\lambda)$, and by considering the pre-image 
$A = [0,\lambda^{1/2}]^2$ of the unit square we see
\[ \ell(\lambda)  = \lambda^{-1/2} \times \area (A) \times  \ell \]
where the $\lambda^{-1/2}$ term is length rescaling.

Similar relations, provable in the same way,  will be stated later (\ref{iota-scale},\ref{p-scale}, \ref{pr-scale}) 
 without repeating the proof.

\paragraph{Aspect 4.  Final definition of a SIRSN.}
Starting from FDDs, 
a conceptual and technical issue is how to continue to understand a SIRSN as a process over the whole continuum.  
As an analogy, for continuous-time stochastic processes one typically seeks some 
sample path regularity property such as {\em c{\`a}dl{\`a}g}.  
So one might seek  some notion of ``regularity" for a realization, 
for instance a.e. continuity of routes $\RR(z_1,z_2)$ as $z_1$ and $z_2$ vary. 
A version of such continuity is proved, under extra assumptions, in section \ref{sec-continuity}. 
But as we next explain, in the present context the assumption $p(1) < \infty$ serves as an alternative regularity 
condition that enables us to study global properties of a SIRSN.

There are several possible real-world measures of ``size" of a road segment, 
quantifying the minor road to major road spectrum -- e.g. number of
lanes; level in a highway classification system; traffic volume.  
What about within our model of a SIRSN?
Recalling the definition (\ref{EE-def}) of $\EE(\lambda, r)$, the limit 
\[ \EE(\infty, r):= \cup_{\lambda <\infty} \EE(\lambda, r) \] 
has (because $\cup_{\lambda <\infty} \Xi(\lambda)$ is dense) the interpretation of 
``the set of path elements $\epsilon$ that are on some route $\RR(z_1,z_2)$ with both $z_1$ and $z_2$ at distance $> r$ from $\epsilon$".   
As shown 
in section \ref{sec-stronger}, assumption (\ref{p1intro}) implies that the edge-intensity $p(r)$ of $\EE(\infty,r)$ is finite and scales as 
$p(r) = p(1)/r$.  
Moreover the random process $\EE(\infty, r)$ is independent of the sampling process $(\Xi(\lambda), 0 < \lambda < \infty)$ and is an intrinsic 
part of the global structure of the SIRSN. 
So if we intuitively interpret $\EE(\infty,r)$  as ``the roads of size $\ge r$", then we have a mathematically convenient notion of ``size of a road segment"  
emerging from our setup without explicit design.  
Intuitively, one could view the limit 
$\EE(\infty, 0+) := \cup_{r>0} \ \EE(\infty,r)$ 
as the continuum network of interest. 
But at a technical level it is not clear what are the properties of a realization  of
$\EE(\infty, 0+)$, and we do not study it in this paper.

\section{The binary hierarchy model}
\label{sec-constr1}
The construction of this model, our basic example of a SIRSN,  occupies all of section \ref{sec-constr1}, 
in several steps.
\begin{itemize}
\item A construction on the integer lattice (sections \ref{sec-1-lattice} - \ref{sec-finesse})
\item Extension to the plane (sections \ref{sec-1-rational} - \ref{sec-1-plane})
\item Further randomization to obtain invariance properties (section \ref{sec-1-complete}).
\end{itemize}

\subsection{Routes on the lattice}
\label{sec-1-lattice}

For an integer $x \neq 0$, 
write $\height(x)$ for the largest $j \in \Ints^+$ such that $2^j$ divides $x$; 
in other words the unique $j$ such that $x = (2k+1)2^j$ for some $k \in \Ints$.
Set $\height (0) = \infty$.  
For later use note that
in one dimension, any integer interval $[m_1,m_2]$ contains a {\em unique} 
integer of maximal height, which we call $\peak[m_1,m_2]$.  
For instance $\peak[67,99] = 96$ and $\peak[34,59] = 48$.

Until section \ref{sec-1-rational} we will work on the integer lattice 
$\Ints^2$, with vertices $z = (x,y)$ whose coordinates have heights $\ge 0$.  
 While we are working on the lattice it is convenient to use $L^1$ distance 
$||z_2 - z_1||_1 := |x_2-x_1| + |y_2-y_1|$.  
Note also that until section \ref{sec-finesse} we work with {\em deterministic} constructions.

Write $L^{(X)}_x$ and $L^{(Y)}_y$ for the lines through  $\{(x,y), y \in \Ints\}$ 
and $\{(x,y), x \in \Ints\}$.   
The height of a line $L^{(X)}_x$is the height of $x$. 

Fix a parameter $ 1/2 < \gamma < 1$.
Associate with lines at height $h$ a cost-per-unit length 
equal to $\gamma^h$. 
Now each path in the lattice has a cost, being the sum of the edge costs.
Visualize a road network in which one can travel along a height-$h$ road at speed $1/\gamma^h$; 
 so the cost equals time taken.

Define the route $\rr(z_1,z_2)$ to be a minimum-cost path between $z_1$ and $z_2$.  
There is a uniqueness issue: for instance, for any minimum-cost path from  $(i,i)$ to $(j,j)$ there is an equal cost path obtained by 
reflection $(x,y) \to (y,x)$.  
However, the estimates from here through section \ref{sec-further} hold when $\rr(z_1,z_2)$ is any choice of minimum-cost path.
We will deal with uniqueness in section \ref{sec-finesse}.

A key point of the construction is that if we scale space by $2$ then the scaled
structure 
on the even lattice $(2\Ints)^2$ agrees with the original substructure 
on the even lattice, up to a constant multiplicative factor in edge-costs, and so 
the route between two even points will be the same whether we work in $\Ints^2$ 
or $(2\Ints)^2$.  So this ``invariance under scaling by $2$" property is built into the model at the start.

The fact that moving along the axes has zero cost may seem worrrying
but actually causes no difficulty 
(we will later apply a random translation, and the original axes do not appear in the final process).   
Note that the cost associated with the line segment from $(2^h,2^h)$ to $(2^h,0)$ is $\gamma^h 2^h$ 
and the constraint $\gamma > 1/2$ 
is needed to make this cost increase with $h$.
Intuitively, if $\gamma$ is near $1$ then the route $\rr(z_1,z_2)$ will stay inside
or near 
the rectangle with opposite corners $z_1, z_2$, whereas 
if $\gamma$ is near $1/2$ then the route may go far away from  
the rectangle to exploit high-speed roads.

For this model we will show a property stronger than (\ref{D1}); the ratio of route-length to distance is uniformly bounded.
\begin{Proposition} 
\label{Pstretch}
There is a constant $K_\gamma < \infty$ such that
\[ \len \  \rr(z_1,z_2) \le K_\gamma ||z_2 - z_1||_1, \quad \forall z_1, z_2 \in \Ints^2 . \]
\end{Proposition}

Some intuition about possible paths in this model is provided by Figure 3 
(the reader should imagine the ratios of longer/shorter edge lengths as
larger than drawn).
We might have a route as shown in the figure, where the two long edges are very fast freeways.
But such a route is not possible if the fast freeways are too far from the  start and destination points.   
The latter assertion will follow from Lemma \ref{Lsquare}.

\setlength{\unitlength}{0.1in}
\begin{picture}(20,21)(-5,-2)
\put(20,15){\line(-1,0){20}}
\put(20,15){\line(0,-1){15}}
\put(15,0){\line(1,0){5}}
\put(15,0){\line(0,1){2}}
\put(15,2){\line(1,0){1}}
\put(16,2){\circle*{0.4}}
\put(0,11){\line(1,0){3}}
\put(0,11){\line(0,1){4}}
\put(3,11){\line(0,1){1}}
\put(3,12){\circle*{0.4}}
\end{picture}

{\bf Figure 3.}  
{\small Routes like this are possible.}

\medskip
\noindent
One might expect some explicit algorithmic description of routes $\rr(z_1,z_2)$ that one can use 
to prove the results in sections \ref{sec-tdm} - \ref{sec-1-plane}, but we have been unable to do so. 
Instead our proofs rely on finding internal structural properties that routes must have.

\subsection{Analysis of routes in the deterministic model}
\label{sec-tdm}
Consider the route from $z_1 = (x_1,y_1)$ to $z_2 = (x_2,y_2)$.
The $x$-values taken on the route form some interval $I_x \supseteq
[\min(x_1,x_2), \max(x_1,x_2)]$, 
and similarly the $y$-values form some interval $I_y$. 
Consider the point
$z^* = (x^*,y^*) = (\peak(I_x),\peak(I_y))$ 
and call this point $\peak^{(2)}\rr(z_1,z_2)$.
The notation reminds us that $\peak^{(2)}\rr(z_1,z_2)$ depends on the route $\rr(z_1,z_2)$, which may not be unique.

\begin{Lemma}
\label{Lpeak}
Consider the route $\rr(z_1,z_2)$ from $z_1$ to $z_2$. \\
(i) The route passes through $z^* = \peak^{(2)}\rr(z_1,z_2)$. \\
(ii) The route meets the line $L^{(X)}_{x^*}$ in either the single point $z^*$ or in one line
segment containing $z^*$ (and similarly for $L^{(Y)}_{y^*}$). \\
(iii) Suppose the route passes 
through a point $(x^*,y)$ (for some $y \neq y^*$) 
and through a point $(x,y^*)$  (for some $x \neq x^*$).
Then the route between those points is the two-segment route via $z^*$.  \\
(iv) Suppose $z^* = z_1$. If 
$z_2$ is in a certain quadrant relative to $z_1$, 
for instance the quadrant $[x_1,\infty) \times [y_1,\infty)$, 
then the route from $z_1$ to $z_2$ stays in that quadrant. 
\end{Lemma}
\proof 
We first prove (iii). 
It is enough to prove that, 
amongst routes between $(x^*,y)$ and  $(x,y^*)$, 
 the two-segment route via $z^*$ is
the {\em unique} minimum-cost route.
In order to get from $(x^*,y)$ to the line $L^{(Y)}_{y^*}$ the route must use at least 
$|y - y^*|$ vertical unit edges; by definition of $x^* = \peak(I_x)$, if these edges are
not precisely
the line segment from $(x^*,y)$ to $z^*$ then the cost of these edges will be
strictly larger; and  similarly for horizontal edges.
This establishes the uniqueness assertion above, and hence (iii).

For (i), if the hypothesis of (iii) fails then the route must go through $z^*$, 
whereas if it holds then the conclusion of (iii) implies the route goes through $z^*$.
For (ii), if false then the the route passes through some two points 
$(x^*,y^\prime)$ and $(x^*,y^{\prime \prime})$ but not the intervening points on
that line.
But (as in the argument for (iii)) the minimum cost path between those two points is
the direct line between them.

Finally, (iv) follows from (ii), because if (iv) fails then the route meets one
boundary of the quadrant
in more than one segment.  
\qed

\begin{Lemma}
\label{Lsquare}
The route $\rr(z_1,z_2)$ from $z_1$ to $z_2$ stays within the square of side $K^\prime_\gamma 
\ ||z_2-z_1||_1$ centered at $z_1$, 
where $K^\prime_\gamma$ depends only on $\gamma$.
\end{Lemma}
\proof 
Choose the integer $h$ such that
\[ 2^{h-1} < ||z_2 - z_1||_1 \le 2^h . \]
As illustrated in Figure 4, there is a square of the form 
$S = [(i-1)2^h, (i+1)2^h] \times [(j-1)2^h,(j+1)2^h]$ 
containing both $z_1$ and $z_2$ 
(note here $i$ and $j$ may be even or odd).  
We may suppose the route does not stay within $S$ (otherwise the result is trivial).
For any point $z$ outside $S$, call the $L^\infty$ distance from $z$ to $S$, 
that is the number $d$ such that $z$ is on the boundary of the concentric square 
$S_d = [(i-1)2^h - d, (i+1)2^h + d] \times [(j-1)2^h - d,(j+1)2^h + d]$, 
the {\em displacement} of $z$.   
Now let $d$ be the maximum displacement along the route $\rr(z_1,z_2)$, and choose 
a point $z^\prime$ along the route with displacement $d$.  
So the route stays within $S_d$.

\setlength{\unitlength}{0.07in}
\begin{picture}(30,28)(-20,-3)
\put(0,0){\line(1,0){10}}
\put(0,10){\line(1,0){10}}
\put(10,0){\line(0,1){10}}
\put(0,0){\line(0,1){10}}
\put(4,7){\circle*{0.8}}
\put(7,4.3){\circle*{0.8}}
\put(1.7,6.6){$z_1$}
\put(6.4,2.4){$z_2$}
\put(5,0){\line(0,1){10}}
\put(0,5){\line(1,0){10}}
\put(-10,20){\line(1,0){30}}
\put(-10,20){\line(0,-1){5}}
\put(20,20){\line(0,-1){5}}
\put(5,15){\vector(0,1){4.6}}
\put(5,15){\vector(0,-1){4.6}}
\put(3,15){$d$}
\put(-1,2.5){\vector(0,1){2.1}}
\put(-1,2.5){\vector(0,-1){2.1}}
\put(-4,2.0){$2^h$}
\put(11,1.5){$S$}
\put(21,16.5){$S_d$}
\put(10,21){$z^\prime$}
\put(10,20){\circle*{0.8}}
\end{picture}

{\bf Figure 4.}  
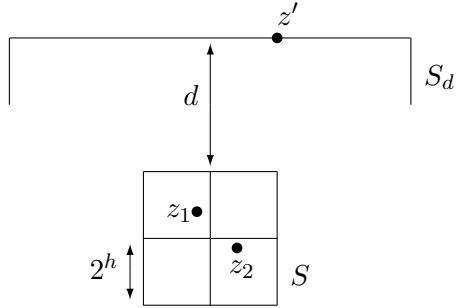
{\small Construction for proof of Lemma \ref{Lsquare}.}

\medskip
\noindent
We may assume, as in Figure 4, that $z^\prime$ is on the top edge of $S_d$.
The route needs to cover the vertical distance $d$ between the top edges of $S$ and $S_d$
{\em twice} (up and down) while staying within $S_d$, 
which has side-length $2^{h+1} + 2d$.   
Now within any integer interval of length $a$ the second-largest height $H$ satisfies 
$2^H \le a$. 
So the cost ($C$,  say) of the route outside $S$ is at least the cost associated with this second-largest height,
which is given by 
\[ d \gamma^H \mbox{ where } 2^H \le 2^{h+1} + 2d.  \] 
Setting $d = b2^h$, 
this inequality implies
\[ \log_2 C \geq \log_2 b + h + (h+1 + \log_2 (1+b)) \log_2 \gamma . \]
But for this to be the minimum-cost path, the cost outside $S$ must be less than the cost of going round the boundary 
of $S$, which is at most $\gamma^h \times 2^{h+2}$.  So
\[ \log_2 C \le h \log_2 \gamma + h + 2 . \]
This inequalities combine to show
\[ \log_2 b + (1 + \log_2 (1+b)) \log_2 \gamma \le 2 \] 
which, because $\gamma > 1/2$,  implies that $b$ is bounded by some constant $b_\gamma$. 
\qed

\medskip 
Lemma \ref{Lsquare} makes Proposition \ref{Pstretch} look very plausible, but to prove it we need to extend Lemma 
\ref{Lpeak} to develop internal structural properties that routes must have.

Call a sequence of integers $i_1, i_2, \ldots, i_m$ 
a {\em height-monotone sequence} from $i_1$ to $i_m$  if \\
(i) $\height(i_1) > \height(i_2) > \ldots > \height(i_m) \ge 0 $;\\
(ii) $|i_{j+1} - i_j| < 2^{\height(i_j)}, \quad 1 \le j < m$.  \\
Suppose, for integers $m_1, m_2, m^*$, we are given  a height-monotone sequence $m^* = i_1, i_2, \ldots, i_m = m_2$ 
and a height-monotone sequence $m^* = j_1, j_2, \ldots, j_q = m_1$. 
Then we can form the concatenation \\
$m_1 = j_q, j_{q-1}, \ldots, j_2, m^*, 1_2, \ldots, i_m = m_2$. 
Call a sequence that arises this way an {\em admissable sequence} from $m_1$ to
$m_2$.  
See Figure 5.

 \begin{center}
 \begin{tabular}{crrrrrrrrr}
 5&&&&&&&96&& \\
 4&&&&&&80&&&\\
 3&&&72&&&&&&\\
 2&&&&&&&&&100\\
 1&&&&74&&&&& \\
 0&&&&&75&&&99&\\
 (height)&&&&&&&&&
 \end{tabular}
 \end{center}

{\bf Figure 5.} 
{\small An admissable path from 75 to 99.  This path has range $100 - 72 = 28$.}

\medskip
\noindent
Regard a height-monotone or admissable sequence as a path of steps where a step from
$i$ to $j$ has 
length $|j-i|$.
It is clear from (ii) that the length of the path in (i) is at most twice 
the length of the first step.
We deduce the following crude bound.
\begin{quote}
(*) The total length of an admissable path is at most 
$4$ times the {\em range} of the path, where the range is the difference between the 
maximum and minimum integer points visited by the path.
\end{quote}
Proposition \ref{Pstretch} follows immediately from Lemma \ref{Lsquare}, the bound (*) above and the 
following lemma.
\begin{Lemma}
\label{Ladmit}
The route from $z_1 = (x_1,y_1)$ to $z_2 = (x_2,y_2)$
consists of alternating horizontal and vertical segments, in which 
the successive distinct $x$-values of the segment ends (the turning points) 
form an admissable sequence from $x_1$ to $x_2$,
and the successive distinct $y$-values form an admissable sequence from $y_1$ to $y_2$.
\end{Lemma}
\proof 
In view of Lemma \ref{Lpeak} we can reduce to the case where 
$z_1 = \peak^{(2)}\rr(z_1,z_2)$, and we need to show that the successive distinct 
$x$-values form a height-monotone sequence, as do the $y$-values.
Without loss of generality suppose that $x_1 \le x_2$, that $y_1 \le y_2$ and that the first segment is horizontal.  
So the route is of the form
\[ (x_1,y_1) = (x_{(1)},y_{(1)}) \to (x_{(2)},y_{(1)}) \to (x_{(2)}, y_{(2)}) \to (x_{(3)}, y_{(2)}) \to \ldots  \]
It suffices to show that for each edge of the route, say the edge 
$(x_{(i)}, y_{(i)}) \to (x_{(i+1)}, y_{(i)})$, and for each point (say $(x^*,y_{(i)})$) on the edge other than the starting point,  we have 
$\height(x^*) < \height (x_{(i)})$.   
This is true for the first two edges of the route by definition of $z_1$ as $ \peak^{(2)}\rr(z_1,z_2)$.  
Suppose it fails first at some point $(x^*,y_{(i)})$.   
Then the route has proceeded 
$(x_{(i)},y_{(i-1)}) \to (x_{(i)},y_{(i)}) \to (x^*,y_{(i)})$ 
instead of the alternate path via $(x^*,y_{(i-1)})$.  
Now inductively $\height(y_{(i)}) < \height(y_{(i-1)})$, 
so the cost of the horizontal edge is less in the alternate path; so for the route to have smaller cost it must happen that 
the cost of its vertical edge is smaller than in the alternate path, that is $\height(x_{(i)}) > \height(x^*)$, 
contradicting the supposed failure. 
\qed

\subsection{Further technical estimates}
\label{sec-further}
The next lemma will be key to bounding network length, 
more specifically to showing $\ell < \infty$ later.
\begin{Lemma}
\label{Lbox}
There exists an integer $b \ge 1$, depending only on $\gamma$,
such that for all $h \ge 0$ and all rectangles of the form 
$[i2^{h+b}, (i+1)2^{h+b}] \times [j2^h, (j+1)2^h]$, 
the route $\rr(z_1,z_2)$ between two points $z_1, z_2 \in \Ints^2$ outside 
(or on the boundary of) the rectangle does not use any horizontal edge
strictly inside the rectangle.
\end{Lemma}
Note there may be routes using a vertical line straight through the rectangle.

\proof
Suppose false; 
then there are two points $z_1, z_2$ on the boundary of the rectangle such that 
the route between them lies strictly within the rectangle 
 and contains a horizontal edge.  
Because the speed on an interior edge is less than the speed on a parallel boundary
edge, 
this cannot happen when $z_1$ and $z_2$ are in the same or adjacent boundaries of
the rectangle, 
because the path around the boundary is faster.
Suppose they are on the top and the bottom boundaries. 
Then the height of the horizontal edge is less than the heights of the starting and 
ending $y$-values, contradicting Lemma \ref{Ladmit}.
The only remaining case is when $z_1$ and $z_2$ are on the left and right boundaries.  
Using Lemma \ref{Ladmit} again, the route cannot use a vertical edge inside the rectangle, 
so the only possibility is a single horizontal segment passing through the rectangle.  
Such a path has cost at least $2^{h+b} \  \gamma^{h-1}$, 
because the height of the line is at most $h-1$,
whereas the path around the boundary has cost at most 
$2^{h+b} \  \gamma^{h} + 2^h \  \gamma^{b+h}$. 
So the potential route is impossible when 
$2^b + \gamma^b < 2^b \gamma^{-1}$ 
which holds for sufficiently large $b$.
\qed

\begin{Corollary}
\label{Crh}
If a route $\rr(z_1,z_2)$ uses a height-$h$ segment through $z_0$, then 
$\min(||z_1-z_0||_1 , ||z_2-z_0||_1) \le 2^h(2^b + 1) $
for $b$ as in Lemma \ref{Lbox}.
\end{Corollary}
\proof 
Consider a unit-length horizontal (without loss of generality) 
edge of height $h$ at $z_0$.  
It is in the interior of some rectangle of the form 
$[i2^{h+1+b}, (i+1)2^{h+1+b}] \times [j2^{h+1}, (j+1)2^{h+1}]$.
By Lemma \ref{Lbox} applied with $h+1$, 
either $z_1$ or $z_2$ must be within that rectangle.
\qed

Perhaps surprisingly, we do not make much explicit use of the 
deterministic function 
$\costbf(z_1,z_2)$ giving the cost of the minimum-cost route in the integer lattice, 
but will need the following bound.
\begin{Lemma}
\label{Lcost}
There exists a constant $K^{\prime \prime}_\gamma$ such that
\[ \costbf(z_1,z_2) \le K^{\prime \prime}_\gamma ||z_2 - z_1||^\beta \]
where $\beta := \log (2 \gamma) / \log 2$.
\end{Lemma}
\proof
As in Figure 4 in the proof of Lemma \ref{Lsquare},
 there is a square of the form 
$S = [(i-1)2^h, (i+1)2^h] \times [(j-1)2^h,(j+1)2^h]$ 
containing both $z_1$ and $z_2$, where $h$ is the  
 integer such that
$ 2^{h-1} < ||z_2 - z_1||_1 \le 2^h $.  
As observed there, the cost of going all around the boundary of $S$ is 
$O(\gamma^h 2^h)$.   
By considering a  path from $z_1$ using the ``greedy" rule of always switching to an orthogonal line of greater height, 
it is easy to check that the cost of this greedy path from $z_1$ to the boundary of $S$ is also $O(\gamma^h 2^h)$.
Hence $\costbf(z_1,z_2) = O(\gamma^h 2^h)$ and the result follows.

\subsection{Finessing uniqueness by secondary randomization}
\label{sec-finesse}
As previously observed, minimum-cost paths are not always unique. 
We conjecture that, at least when $\gamma$ is not algebraic, there is some 
simple classification 
of when and how non-uniqueness occurs. 
But instead of addressing that issue we can finesse it by introducing randomness (which we need later, anyway) at this stage.  
One possible way to do so would be to use the uniform distribution on minimum-cost paths.  
Instead we use what we will call 
{\em secondary randomization} to choose between non-unique minimum-cost paths. 
Place i.i.d. Normal$(0,1)$ random variables (``weights") $\zeta_e$ on the edges $e$
of $\Ints^2$.  Any path has a weight 
$\sum_{e \mbox{ in path}} \zeta_e$.  
Define the route $\RR_0(z_1,z_2)$ to be the minimum-weight path in the set of
minimum-cost paths from $z_1$ to $z_2$.

\subsection{Extension to the binary rational lattice}
\label{sec-1-rational}
The notion of {\em height} extends to binary rationals: 
if $x \in \Reals$ is a binary rational and $x \neq 0$, 
write $\height(x)$ for the largest $j \in \Ints$ such that $2^j$ divides $x$; 
in other words the unique $j$ such that $x = (2k+1)2^j$ for some $k \in \Ints$.

For $- \infty < H <  \infty$ let $\Ints^2_{H}$ be the lattice on vertex-set 
$\{2^H z: \  z \in \Ints^2\}$, in other words on the set of points in $\Reals^2$ whose coordinates have height $\ge H$.
So far we have been working on the integer lattice $\Ints^2$, but now the results we have proved extend by (binary) scaling to 
analogous results on the lattices $\Ints^2_{H}$.  
We will use such scaled results as needed.

Note in particular the following consistency condition as $H$ varies. 
Take $H_1 < H_2$.  
Consider the route, in $\Ints^2_{H_1}$, between two vertices of $\Ints^2_{H_2}$.  
By Lemma \ref{Ladmit} and the definition of admissable, any minimum-cost path stays within the lattice $\Ints^2_{H_2}$.  
So the set of minimum-cost paths is the same whether we work  in $\Ints^2_{H_1}$ or in $\Ints^2_{H_2}$. 
Note also that each edge $e$ in $\Ints^2_H$ corresponds to two edges 
$e_1, e_2$ of $\Ints^2_{H-1}$.  
So we can couple the edge-weights by making $\zeta_e = \zeta_{e_1} + \zeta_{e_2}$ 
(only this infinite divisibility property of the Normal is relevant to the
construction) and this gives a 
``consistency of secondary weights" property, which implies that 
the random route  $\RR_0(z_1,z_2)$
is also the same whether we work  in $\Ints^2_{H_1}$ or in $\Ints^2_{H_2}$.

So we have now defined random routes $\RR_0(z_1,z_2)$ for all unordered pairs of vertices 
 in $\Ints^2_{-\infty}:= \cup_{H>-\infty} \Ints^2_{H}$. 
From the ``minimality" in the construction it is clear that the routes
satisfy the route-compatability properties (iii,iv) from section  \ref{sec-TS}.

\subsection{Extension to the plane}
\label{sec-1-plane}

We want to define routes $\RR_0(z_1,z_2)$ between general points $z_1,z_2$ of
$\Reals^2$ as $H \to - \infty$ limits of the routes 
$\RR_0(z_1^H,z_2^H)$ between vertices such that 
\begin{equation}
z_i^H \in  \Ints^2_H, 
\ \  z_i^H  \to z_i \ \ (i = 1,2)
\label{ziH}
\end{equation}
Proposition \ref{PnullA} formalizes this idea.  
The proof in this section  is the most intricate part of the construction, which can thereafter be completed (section \ref{sec-1-complete}) by ``soft" arguments.

As a first issue, 
what does it mean to say that, under (\ref{ziH}), 
\begin{equation}
\mbox{routes
$\rr(z_1^H,z_2^H)$ converge to a route $\rr(z_1,z_2)$?}
\label{routes-H}
\end{equation} 
We define this to mean:
\begin{quote}
for each $H_0 > - \infty$, the subroute $\rr_{H_0}(z_1^H,z_2^H)$
consisting of path segments of $\rr(z_1^H,z_2^H)$ within lines of height $\ge H_0$ is, 
for sufficiently large negative $H$, a path not depending on $H$ -- call this path 
$\rr_{H_0}(z_1,z_2)$.  
\end{quote}
When this property holds, Lemma \ref{Ladmit} implies that $\rr_{H_0}(z_1,z_2)$ is a connected path,
consistent as $H_0$ decreases, and using Proposition \ref{Pstretch} and scaling we
see that the closure of 
$\cup_{H_0 > - \infty} \rr_{H_0}(z_1,z_2)$ defines a route $\rr(z_1,z_2)$ satisfying the  ``jagged"
condition of section \ref{sec-TS}.

\begin{Proposition}
\label{PnullA}
There exists a subset $A \subset \Reals^2$ of zero area such that,  if $z_1$ and
$z_2$ are outside $A$, there exists a random jagged route  $\RR_0(z_1,z_2)$ such that,
whenever (\ref{ziH}) holds, then
(\ref{routes-H}) holds.
\end{Proposition}
The proof relies on the fact that, for particular configurations
illustrated in Figure 6, routes from a certain neighborhood to distant
destinations must
 pass through a particular point.  
 In fact all that matters is the {\em existence} of such a configuration, 
not the particular one we now exhibit. 
Consider a square
$G = [2^h ,2^h +2]^2$ 
and points $b = (2^h +1 , 2^h)$ and  $d = (2^h + 1,2^h + 2^{-h})$,
illustrated in Figure 6.  
What is relevant is the heights of the lines involved,
 indicated in the figure.

\setlength{\unitlength}{0.15in}
\begin{picture}(18,21)(-5,-3)
\put(-1,0){\line(1,0){17}}
\put(-1,2){\line(1,0){1}}
\put(-1,8){\line(1,0){17}}
\put(-1,16){\line(1,0){17}}
\put(0,-1){\line(0,1){17}}
\put(8,-1){\line(0,1){17}}
\put(16,-1){\line(0,1){17}}
\put(7,2){\line(1,0){1}}
\put(7,2){\line(0,1){1}}
\put(7,3){\line(1,0){1}}
\put(-2,-0.3){h}
\put(-2,15.7){1}
\put(-2.15,1.7){-h}
\put(-2,7.8){0}
\put(-0.2,-2){h}
\put(7.8,-2){0}
\put(15.8,-2){1}
\put(7.2,2.13){$\Sigma$}
\put(7.3,7){{\huge G}}
\put(8,8){\circle*{0.3}}
\put(8.2,8.3){{\footnotesize O}}
\put(8,0){\circle*{0.3}}
\put(8.2,0.2){b}
\put(8,2){\circle*{0.3}}
\put(8.23,2.0){d}
\put(6.6,2){\vector(0,1){1}}
\put(5.9,2.25){$\eps$}
\put(15.2,16){\circle*{0.3}}
\put(14.9,16.3){$c^*$}
\put(8,16){\circle*{0.3}}
\put(7.8,16.4){$b^*$}
\put(16,16){\circle*{0.3}}
\put(16.3,15.8){{\footnotesize NE}}
\put(16,0){\circle*{0.3}}
\put(16.3,-0.2){{\footnotesize SE}}
\put(0,16){\circle*{0.3}}
\put(-0.6,16.3){{\footnotesize NW}}
\put(0,0){\circle*{0.3}}
\put(0.2,0.2){{\footnotesize SW}}
\put(16,8){\circle*{0.3}}
\put(16.3,8){{\footnotesize E}}
\put(0,8){\circle*{0.3}}
\put(0.3,8.2){{\footnotesize W}}
\end{picture}

{\bf Figure 6.}
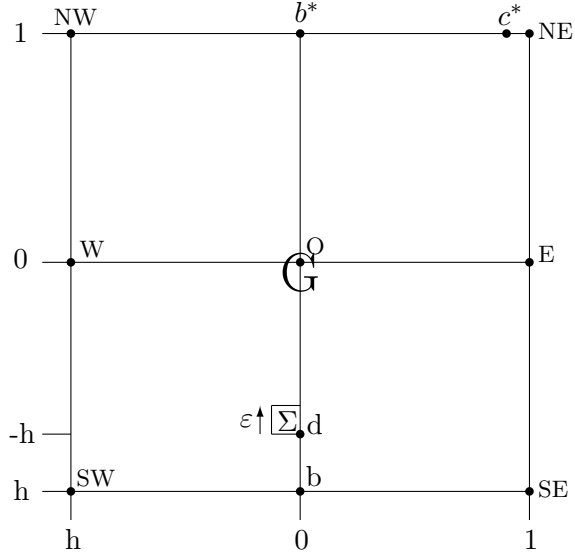
{\small  The big square $G$ and the small square $\Sigma$.  
Marginal labels attached to lines are line-heights, not coordinates.}

\begin{Lemma}
\label{L-key} 
There exist large $h$ and small $\eps$ (depending on $\gamma$) such that, 
in the configuration shown in Figure 6, every route from inside the small square $\Sigma:= d + [ - \eps,0] \times [0,\eps]$  to the
boundary of $G$ passes via $b$.
\end{Lemma}
\proof
For each point $c$ on the boundary of $G$ there is a 
counter-clockwise path $\pi_1(b,c)$ and a clockwise path $\pi_2(b,c)$ 
along the boundary from $b$ to $c$.  
These paths have equal cost for the point 
$c^* = (2^h + 2 - \gamma^{h-1}, 2^h + 2)$, 
which is near the NE corner point of $G$.  
We will need the following lemma.
\begin{Lemma}
\label{Lccc}
There exists $h$ such that the following hold.\\
(a) The paths $\pi_1(b,c^*)$ and $\pi_2(b,c^*)$ attain the minimum cost 
over all paths from $b$ to $c^*$, and are the only paths to do so. \\
(b) The only minimum-cost paths from $d$ to $c^*$ are the two paths consisting of 
the segment $[d,b]$ and the paths $\pi_1(b,c^*)$ or $\pi_2(b,c^*)$. \\
(c) There exists $\eta > 0$ such that any path from $d$ to $c^*$ that avoids the 
segment $[d,b]$ has cost at least $\eta$ greater that the minimum-cost paths.
\end{Lemma}
Note a technical point.  We are working on $\Ints^2_{-\infty}:= \cup_{H>-\infty} \Ints^2_{H}$ 
and $c^*$ may not be in $\Ints^2_{-\infty}$.  To be precise we should replace $c^*$ in the arguments below 
by a sequence $c^*_H \to c^*$,  but  that requires awkward notation we prefer to avoid.

Granted Lemma \ref{Lccc} we deduce Lemma \ref{L-key}  as follows.
Consider a point $c$ on the counter-clockwise path from $b$ to $c^*$ 
(the clockwise case is similar). 
Then the following must hold, because any counter-example path to 
$c$ could be extended along the boundary from $c$ to $c^*$ and would give a 
counter-example to Lemma \ref{Lccc}. \\
(a) The path $\pi_1(b,c)$ is the unique minimum-cost path from $b$ to $c$.\\
(b) The path consisting of the segment $[d,b]$ 
and the path $\pi_1(b,c)$ is the unique minimum-cost path from $d$ to $c$.\\
(c) Any path from $d$ to $c$ that avoids the 
segment $[d,b]$ has cost at least $\eta$ greater that the minimum-cost path. \\
Lemma \ref{Lcost} extends by scaling to 
$\cup_{H>-\infty} \Ints^2_{H}$, and so the function 
$\costbf(\cdot,\cdot)$ extends to a continuous function on $\Reals^2$.
So we can choose $H$ 
so that the square $\Sigma = d + [-2^{-H},0]\times[0,2^{-H}]$  
satisfies  
$\sup_{s \in \Sigma} \costbf(s,d) \le \eta/3$.

It is easy to check that a minimum-cost path from $s \in \Sigma$ to $d$ does not meet
$[d,b]$ 
except at $d$.

Consider $s \in \Sigma$ and a point $c$ as above.
So $\costbf(s,c) \le \costbf(d,c) + \eta/3$.
Suppose a minimum-cost path from $s$ to $c$ does not
meet the segment $[d,b]$.  
Then the path from $d$ to $c$ via $s$ would have cost 
$\le \costbf(d,c) + 2\eta/3$ and would not meet $[d,b]$,
contradicting (c).
So a minimum-cost path from $s$ to $c$ 
must meet the segment $[d,b]$.  
Then, by uniqueness in (b) (for the path from $d$ to $c$),  it must continue via $b$, establishing Lemma \ref{L-key}. 
\qed

\medskip

\noindent
{\bf Proof of Lemma \ref{Lccc}.} 
(I thank Justin Salez for completing the details of this proof.)
In outline, we use the ``structure of paths" results in Lemmas \ref{Lpeak} and \ref{Ladmit}
to reduce to comparing costs of a finite number of 
possible routes. We will make use of the following preliminary observations, which are straightforward to check :\\
(i) the only minimum-cost path from $SE$ to $NW$ is $SE\to SW \to NW$ ; \\
(ii) the only minimum-cost paths from $SW$ to $NE$ are $SW\to NW \to NE$ and $SW\to SE \to NE$ ;\\
(iii) $O\to b^* \to NE$ is a minimum-cost path from $O$ to $NE$.

\medskip
\noindent
Consider assertion (a).  The cost associated with paths $\pi_1(b,c^*)$ and $\pi_2(b,c^*)$ equals $2\gamma^h+2\gamma$, which (by choosing $h$ large) is less than $2$. Now consider some minimum-cost path $\pi$ from $b$ to $c^*$. Since both end-points have their $y-$coordinate at height $\geq 1$, all horizontal segments of $\pi$ must have height $\geq 1$ (Lemma \ref{Ladmit}). In other words, the length of every vertical segment must be an even integer. If the last vertical segment of $\pi$ were ending strictly between $NW$ and $NE$, then its cost would be at least $2$, contradicting optimality. Thus, $\pi$ must pass through $NW$ or $NE$, 
 and observationss (i) or (ii) complete the proof. 

\medskip \noindent
Now consider assertions
(b) and (c).  Let $\pi$ be any path from $d$ to $c^*$, and let $z$ be the point at which $\pi$ first meets the boundary of the rectangle formed by $\{E,W,SW,SE\}$. Let $\pi',\pi''$ denote the subpaths of $\pi$ from $d$ to $z$ and from $z$ to $c^*$, respectively. There are four possible cases :
\begin{itemize}
\item $z\in (E,W)$ : since all segments in $\pi'$ have height $\leq 0$, replacing $\pi'$ by $d\to O\to z$ cannot increase the overall cost. In the resulting path, one may further replace the subpath from $O$ to $c^*$ by $O\to b^* \to c^*$ without increasing the cost, by (iii). This shows :
$$\cost(\pi) \geq 2+\gamma-2^{-h}-\gamma^h.$$
\item $z\in (W,SW)$ : 
all horizontal segments in $\pi'$ have height $\leq -1$, so $\cost(\pi')\geq \gamma^{-1}.$ By (ii), one also has $\cost(\pi'')\geq \cost(z\to NW\to c^*).$ Combining these two facts yields
$$\cost(\pi) \geq 2\gamma+\gamma^{-1}.$$
\item $z\in(SE,E)$ : replacing the subpath $\pi''$ by $z\to NE \to c^*$ cannot increase the overall cost, by part (a). In the resulting path, the subpath from $d$ to $E$ costs at least $\gamma^{-1}+\gamma(1-2^{-h})$, because the horizontal and vertical heights are $\leq -1$ and $\leq 1$, respectively. Thus, 
$$\cost(\pi) \geq 2\gamma+\gamma^{-1}.$$
\item $z\in (SW,SE)$ : all segments of $\pi'$ have height $\leq-1$ except those included in $[O,b]$, which have height $0$. Thus,  
$$\cost(\pi')-\cost(d\to b\to z) \geq (\gamma^{-1}-1)\len\left([b,d]\setminus\pi'\right).$$
Moreover, by part (a), 
$\cost(b\to z)+\cost(\pi'')\geq \cost(\pi_1).$
Thus,
$$\cost(\pi)\geq (\gamma^{-1}-1)\len\left( [b,d]\setminus\pi\right)+\cost(d\to b)+\cost(\pi_1).$$
\end{itemize}
Let us sum up: in the first three cases, the cost of $\pi$ exceeds that of our two candidates by at least $1$, for $h$ sufficiently large. In the fourth case, the excess is at least $(\gamma^{-1}-1)\len\left( [b,d]\setminus\pi\right)$. This proves both (b) and (c), with $\eta = 2^{-h}(\gamma^{-1}-1)$.  
\qed

\medskip
\noindent
{\bf Proof of Proposition \ref{PnullA}.}
In each basic
$2^{h+1} \times 2^{h+1}$ square $G$ of $\Ints^2_{h+1}$ there is a copy of the Figure 6 configuration; let $\Sigma_G$ be 
the corresponding small square. 
Let $B: = \cup_G \Sigma_G$ be the union of those squares, for the fixed $h$ given by  Lemma \ref{L-key}.
Then for $i \geq 1$ let 
$B_i:= \sigma_{2^{-i}}B$ be rescalings of $B$.  Each $B_i$ has the same density, 
which by Lemma \ref{L-key} is non-zero, 
and a straightforward use of the second Borel-Cantelli lemma  
(with sufficiently well-spaced values of $i$) shows that the set 
\[ A:= \{z \in \Reals^2:  \ z \mbox{ in only finitely many } B_i\} \] has area zero.

Now consider $z_1 \in A^c$.
Then there exists a sequence  $i_j = i_j(z_1) \to \infty$ such that $z_1 \in B_{i_j}$ and the associated
$b_{i_j}(z_1) \to z_1$.  
Consider $z_2 \neq z_1$ and 
$(z^H_1,z^H_2) \to (z_1,z_2)$ as in (\ref{ziH}).  
For $j$ larger than some $j_0(z_1,z_2)$, Lemma \ref{L-key} implies that for all sufficiently large $H$ the route 
$\RR_0(z_1^H,z_2^H)$ passes through $b_{i_j}(z_1)$. 
But the routes between the $b_{i_j}(z_1), j \ge 1$ are specified by the construction on $\cup_{H>-\infty} \Ints^2_{H}$.
It follows that,  when $z_1$ and $z_2$ are both in $A^c$ we have convergence in the sense 
of (\ref{routes-H}) to a route $\RR_0(z_1,z_2)$.
\qed

\medskip

We digress to give the technical estimate that will show $\ell < \infty$ in this model.
\begin{Lemma}
\label{LS*}
For the routes $\RR_0$ in Proposition \ref{PnullA}, 
take the union over points $\xi, \xi^\prime$ of a rate-$1$ Poisson point process $\Xi(1)$ of the routes $\RR_0(\xi,\xi^\prime)$, and let 
$\SS^*$ be the intersection of that union with the interior of a unit square 
$U = [i,i+1] \times [j,j+1]$.  
Then the expected length of $\SS^*$ is at most $2^{b+2}$, for $b$ as in  Lemma \ref{Lbox}.
\end{Lemma}
\proof
Lemma \ref{Lbox} was stated for $h \ge 0$ and vertices in $\Ints^2$, but by scaling it holds for $h < 0$ and vertices in $\Reals^2$.
Consider $h < 0$.
Within $U$ there are $2^{-h-1}$ horizontal unit-length line segments  at height $h$, and these can be split into $2^{-2h-1}$ segments of length $2^h$. 
Consider such a line segment, $\zeta$ say.  
It is in the interior of some rectangle of the form 
$[i2^{h+1+b}, (i+1)2^{h+1+b}] \times [j2^{h+1}, (j+1)2^{h+1}]$.
By Lemma \ref{Lbox} applied with $h+1$, the only possible way that the segment $\zeta$ can be in a 
route $\RR_0(\xi,\xi^\prime)$ is if 
$\xi$ or $\xi^\prime$ is within the rectangle.
(And the same holds for any piece of $\zeta$, by considering a sub-rectangle).
The chance the Poisson process contains such a point is at most the area  of the rectangle, which is $2^{2h+2+b}$. 

So the contribution to mean length from a particular segment $\zeta$  is at most $2^h \times 2^{2h+2+b}$, and then the contribution from
height-$h$ horizontal lines is at most
$2^h \times 2^{2h+2+b} \times 2^{-2h-1} = 2^{h+1+b}$.
Summing over $h \le -1$ and adding the same contribution from vertical lines gives the bound  $2^{2+b}$.
\qed

\subsection{Completing the  construction by forcing invariance}
\label{sec-1-complete}
Proposition \ref{PnullA} gives paths $\RR_0(z_1,z_2)$ when $z_1, z_2 \in A^c$.
The process $\RR_0$ cannot be translation- or rotation-invariant (in distribution), because the axes play a special role (infinite speed); though by construction 
the process is invariant under $\sigma_2$ (scaling space by a factor  $2$).
But there is a standard way of trying to make translation-invariant random processes 
out of deterministic processes, by taking weak limits of random translations of the original process.
In our setting this can be done fairly explicitly as follows.  
For $u \in \Reals^2$ let $T_u$ be the translation map 
$T_u(z) = u + z, \ z \in \Reals^2$ on points, and let $T_u$ act on routes in the natural way.
Take $U_n$ uniform on the square $[0,2^n]^2$, and couple the random variables
$(U_n, n \geq 1)$ by setting 
$U_n = U_{n+1} \bmod 2^n$ coordinatewise. 
Define
\begin{equation}
 \RR^{(n)}(z_1,z_2) = T_{- U_n} (\RR_0(z_1+U_n,z_2+U_n)) . 
 \label{Rn-def}
 \end{equation}
In words, translate points  by $U_n$, use  $\RR_0$ to define a route between the translated points, 
and then translate back to obtain a route between
the original points.

Now the only way that $\RR^{(n+1)}(z_1,z_2)$ could be different from 
 $\RR^{(n)}(z_1,z_2)$ is if the route 
$\RR_0(z_1 + U_n, z_2+U_n)$ 
intersects the boundary of the square  $[0,2^n]^2$, 
which, using Lemma \ref{Lsquare}, has chance $O(2^{-n})$.
So we can define a random network $\RRti$ via the a.s. limits
\begin{equation}
\RRti(z_1,z_2) = \RR^{(n)}(z_1,z_2) \mbox{ for all sufficiently large } n . 
\end{equation}
This process is translation-invariant, because for fixed $z \in \Reals^2$ 
the variation distance between the distributions of 
$U_n$ and $U_n + z \bmod 2^n$ tends to zero.

For $0<c<\infty$ write $\sigma_c$ for the scaling map $z \to cz$ on $\Reals^2$, and
recall that 
$\RR_0$ 
is invariant under $\sigma_{2}$. 
Now for $\RR^{(n)}$ at (\ref{Rn-def}),
\begin{eqnarray*}
\sigma_2 \RR^{(n)}(z_1,z_2) 
&=& \sigma_2 T_{-U_n} \RR_*(z_1 + U_n, z_2 + U_n) \\
&=& T_{-2U_n} \sigma_2 \RR_*(z_1 + U_n, z_2 + U_n) \\
&\ed & T_{-2U_n} \RR_*(2z_1 + 2U_n, 2z_2 + 2U_n) \mbox{ by  invariance of $\RR_0$ under $\sigma_2$}\\
&\ed & T_{-U_{n+1}} \RR_*(2z_1 + U_{n+1}, 2z_2 + U_{n+1}) \mbox{ because } U_{n+1}
\ed 2U_n \\
&=& \RR^{(n+1)}(2z_1,2z_2). 
\end{eqnarray*}
Hence the distribution of the limit $\RRti$ is invariant under $\sigma_2$.

Of course our construction so far is not rotationally invariant, but applying a 
uniform random rotation to $\RRti$ gives a network $\RRri$ whose
distribution is 
invariant under rotation, as well as preserving distributional invariance 
under translation and under $\sigma_2$.  
Finally, we get a process $\RR$ with scale-invariant distribution by random 
rescaling via the scale-free distribution: 
\begin{equation}
 \RR(z_1,z_2)  = \sigma_{1/C} \RRri(Cz_1,Cz_2) , \quad 
\Pr (C \in dc) = \sfrac{1}{c \ \log 2}, \ 1 < c < 2 . 
\label{scale-free}
\end{equation}

This completes the construction of  the {\em binary hierarchy model} $\RR$.
To check it satisfies the formal setup of a SIRSN in section \ref{sec-TS}, 
the only remaining issue is to check that the parameters 
$\Ex D_1, \ell$ and $p(1)$ are finite.  
For the former, Proposition \ref{Pstretch} implies the corresponding bound in terms of Euclidean distance
\[ \len \  \RR_0(z_1,z_2) \le 2^{1/2}  K_\gamma ||z_2 - z_1||_2 \] 
and this bound is unaffected by the transformations taking $\RR_0$ to $\RR$.
So $\Ex D_1 \le 2^{1/2}  K_\gamma$.
For $\ell$, in the notation of Lemma \ref{LS*}, 
the edge-intensity of $\cup_{\xi, \xi^\prime \in \Xi(1)} \RR_0(\xi, \xi^\prime)$ is at most 
$2^{b+2} + 2$, the ``$+2$" terms arising from the edges of $\Ints^2$.  
This edge-intensity is unaffected by the transformations taking $\RR_0$ to $\RRri$. 
Scaling by $C$ in (\ref{scale-free}) 
multiplies edge-intensity by $C$, so finally 
$\ell \le (2^{b+2} + 2) \Ex C$.  
To bound $p(1)$, set 
$r(h) = 2^h(2^b +1)$.
 Corollary \ref{Crh} implies that, for routes $\RR_0$, if an edge element 
is in a route between some two points at distance $\ge r(h)$ from the element, 
then the edge has height $\ge h$.  
The edge-intensity of edges with height $\ge h$ equals $2^{1-h}$.  
These quantities are unaffected by translation and rotation; and the scaling 
by $\sigma_C$ can at most increase the edge-intensity by $4$.
So the edge-intensity in $\RR$ of $\EE(\lambda, r(h))$ is 
$p(\lambda, r(h)) \le 4 \cdot 2^{1-h}$ 
Choosing $h$ such that $r(h) < 1$ we deduce $p(1) < \infty$.

\subsection{Remarks on section \ref{sec-constr1}.} 
The ``combinatorial" arguments in sections \ref{sec-1-lattice} - \ref{sec-finesse}
are obviously specific to this model.  But the property implicit in 
Lemma \ref{L-key} 
(that there exist configurations in which all long routes from a small 
neighborhood exit the unit disc at the same point) 
is closely related to desirable
structural properties of SIRSNs 
discussed in section \ref{sec-cul}. 

Lemma \ref{Lellp} later shows that in general $\ell \le 2p(1)$,
 so our argument above that $\ell < \infty$ could be omitted, 
though it is pleasant to have a self-contained construction.

\section{Other possible constructions}
\label{sec:otherC}
The model in section 3 has some very special features, in particular that in any realization we see a (scaled and rotated) square lattice of roads.
Below we outline two other constructions which,  we conjecture, produce SIRSNs,   the technical dificulty being to prove a.s. uniqueness of routes defined as minimum-cost paths.

\subsection{The Poisson line process model}
\label{sec-constr2}
For each $m = 1,2,3,\ldots$ take a rate-$1$ Poisson line process, and attach
Uniform$(m-1,m)$ 
marks to the lines; the union of all these  
is a Poisson line process with ``mark measure" being Lebesgue measure on
$(0,\infty)$. 
By a one-to-one mapping of marks one can transform to the mark measure with density 
$x^{- \gamma}$ on $0<x<\infty$, 
where we take the parameter $2 < \gamma < \infty$.  
So in any small disc, there is some finite largest mark amongst lines intersecting
the disc.

Picturing the lines as freeways and the marks as speeds, for any pair of points
$z_1, z_2$ on the lines there is some finite minimum time 
$\time(z_1,z_2)$ over all routes from $z_1$ to $z_2$, and 
analogous to Lemma \ref{Lcost} 
one can show 
(Wilf Kendall: personal communication) that this function extends to a random continuous function $\time(z_1,z_2)$ on the plane. 
The technical difficulty is to show that for given $(z_1,z_2)$ there is an a.s. unique route attaining that time; if that were proved, establishing the remaining properties required of an SIRSN would be straightforward.  
In particular, scale-invariance would follow from the form $x^{-\gamma}$ of the mark density.

\subsection{A dynamic proximity graph model}
\label{sec-constr3}
This potential construction of a SIRSN is based on 
a space-time Poisson point process $(\Xi(\lambda), 0 < \lambda < \infty)$. 
Note that to study such a SIRSN one would use an independent 
Poisson point process to define
$\SS(\lambda)$.
Note also that the corresponding ``static" model,  
called the {\em Gabriel} network, is a member of the family of
{\em proximity graphs} described in 
\cite{jaromczyk,me-spatial-1};  any family member 
could be used in the construction below.

Here's the construction rule.
\begin{quote}
When a point $\xi$ arrives at time $\lambda$, consider in turn each existing 
point $\xi^\prime \in \Xi(\lambda-)$, and create an edge 
$(\xi,\xi^\prime)$ if the disc with diameter $(\xi,\xi^\prime)$ contains no other 
point of $\Xi(\lambda-)$.
\end{quote}
Write $\GG(\lambda)$ for the time-$\lambda$ network on points $\Xi(\lambda)$.  
Note the automatic scale-invariance property
\begin{quote}
the action of $\sigma_c$ on $\GG(\lambda)$ gives a network distributed as 
$\GG(c^{-2}\lambda)$.
\end{quote}
Now fix a parameter $0 \le \gamma < \gamma_*$ for some sufficiently 
small $\gamma_* > 0$ and view an edge created at time $\lambda$ as a road with speed 
$\lambda^{-\gamma}$.
Defining routes in $\GG(\lambda)$ as minimum-time paths, it seems intuitively 
plausible, 
as in the Poisson line process model, that that we can extend the minimum-time function on $\cup_\lambda \Xi(\lambda)$ to a continuous function $\time(z_1,z_2)$ and then prove  there is an a.s. unique route attaining that time.  Again, if that were proved, establishing the remaining properties required of an SIRSN would be straightforward.  
In particular, scale-invariance would follow from the fact that the 
construction rule is scale-invariant.

 \section{Properties of weak SIRSNs}
 \label{sec-properties}
 In this section we study properties that hold for any weak SIRSN, 
 that is when we do not require (\ref{p1intro})  but instead require (\ref{ell-finite}). 
 These are essentially properties of the sampled subnetworks $\SS(\lambda)$ for fixed $\lambda$ -- we cannot get 
 $\lambda \to \infty$ results.

 \subsection{No straight edges at typical points}
 \label{sec-jagged} 
 If a point $\xi$ of $\Xi(\lambda)$ is the start of some straight line segment of length $\ge r$ in $\SS(\lambda)$ 
 then consider the subroutes of length exactly $r$ from $\xi$.  
The edge process of such subroutes has some
edge-intensity $\iota(\lambda,r)$.
In independent copies of $\Xi(1)$ these edge-processes cannot have any positive-length overlap. 
So by regarding $\Xi(n)$ as the union of  $n$ independent copies of $\Xi(1)$ we 
have $\iota(n,r) = n \iota(1,r)$.
But by the general scaling property (\ref{SS-def-4})
\begin{equation}
\iota(\lambda,r) = \lambda^{1/2} \iota(1,r \lambda^{1/2}). 
\label{iota-scale}
\end{equation} 
Since $\iota(1,r) \le \ell < \infty$ these two different scaling relations imply 
$\iota(1,r) = 0$ for all $r>0$.

 This proves (a) below; note the consequence (b), implied by the definition of {\em feasible path} in the 
 section \ref{sec-TS} setup.
 \begin{Proposition} 
 \label{Pjagged}
 $S(\lambda)$ has the following properties a.s. \\
 (a) $S(\lambda)$ contains no line segment $[\xi,z]$ of positive length, for any $\xi \in \Xi(\lambda)$. \\
 (b) The route $\RR(\xi_1,\xi_2)$ between two points of $\Xi(\lambda)$ does not pass through
 any third point $\xi_3$ of $\Xi(\lambda)$.
 \end{Proposition}
 
 \subsection{Singly and doubly infinite geodesics}
 Recall from section \ref{sec-TS-2} that
 a {\em singly infinite geodesic} from a point $\xi_0$ in $\SS(\lambda)$ is an infinite path, starting from $\xi_0$, 
 such that any finite portion of the path is a subroute of some route $\RR(\xi_0,\xi)$.
 Lemma \ref{Lss1} showed
  \begin{equation}
 \label{Lsingly}
 \mbox{There is a.s. at least one singly infinite geodesic from each point of $\SS(\lambda)$.}
 \end{equation}
A {\em doubly infinite geodesic} in $\SS(\lambda)$ is a path $\pi$ which
 is an increasing union of segments $\pi_k$, where each $\pi_k$ is a segment of
some route 
$\RR(\xi_k,\xi^\prime_k)$ between two points of $\Xi(\lambda)$, and both endpoints of $\pi_k$ go to infinity.

Previous work on very different (e.g. percolation-type \cite{MR1387641}) networks
suggests there may be a general principle:
\begin{quote}
In natural models of random networks on $\Reals^2$ or $\Ints^2$, 
doubly infinite geodesics do not exist.
\end{quote} 
Proposition \ref{P-geod} proves this for weak SIRSNs based on a 
simple scaling argument.  Note however this argument depends implicitly
upon our assumption $\ell < \infty$ which seems rather special to our setting.

Recall the setup of (\ref{EE-def}, \ref{plar-def}).
$\EE(\lambda,r) \subset \SS(\lambda)$ is the set of points $z$ in edges of $\SS(\lambda)$ such that 
$z$ is in the route $\RR(\xi,\xi^\prime)$ for some 
 $\xi, \xi^\prime$ of $\Xi(\lambda)$ such that $\min(|z-\xi|, |z-\xi^\prime|) \ge r$.
And  $p(\lambda,r)$ is the edge-intensity of $\EE(\lambda,r)$.
By scaling, 
\begin{equation}
p(\lambda,r) = \lambda^{1/2} p(1,r\lambda^{1/2})
\label{p-scale}
\end{equation}

\begin{Proposition}
\label{P-geod}
$p(\lambda,r) \to 0$ as $r \to \infty$.  
In particular, $\SS(\lambda)$ has a.s. no  doubly infinite geodesics.
\end{Proposition}
\proof 
For fixed $\lambda$ the edge-processes $\EE(\lambda, r)$ can only decrease as $r$ increases, 
and the limit $\EE(\lambda, \infty) := \cap_r \EE(\lambda, r)$ is by definition the set of path elements in 
doubly infinite geodesics.  This limit has edge-intensity
$p(\lambda, \infty) = \lim_{r \to \infty} p(\lambda,r) \ge  0$.  
So it is enough to prove $p(\lambda,\infty) = 0$.  Suppose not. 
Then by the scaling relation (\ref{p-scale})
\[ p(\lambda,\infty) = \lambda^{1/2} p(1,\infty), \quad 0 < \lambda < \infty . \]
We claim that in fact 
\[  \mbox{$\EE(\lambda, \infty) = \EE(1, \infty)$ a.s. for $\lambda < 1$,} 
\] 
which (because we know $p(1,\infty) < \infty$) implies $p(1,\infty) = 0$ and completes the proof.

To prove the claim, note that
for any finite-length segment $\pi_0$ of a doubly infinite geodesic in $\SS(1)$,
there are an infinite number 
of distinct pairs $\xi_j, \xi^\prime_j$ of $\Xi(1)$ such that $\RR(\xi_j,\xi^\prime_j)$ contains
$\pi_0$, and for each pair 
there is chance $\lambda^2$ that both points are in $\Xi(\lambda)$.
 These events are independent (because $\Xi(\lambda)$ is obtained from $\Xi(1)$ by independent sampling) 
so a.s. an infinite number of pairs $\xi_j, \xi^\prime_j$ are in $\Xi(\lambda)$,  implying that $\pi_0$ is in a doubly infinite geodesic
of $\SS(\lambda)$.
\qed

{\bf Remark.}  
The limit used here is different from the limit 
$p(r) := \lim_{\lambda \to \infty} p(\lambda,r)$ 
featuring in assumption (\ref{p1intro}).

\subsection{Marginal interpretation of $\ell$}
\label{ell-formula}
Recall $\ell$ is defined as the edge-intensity of
$\SS(1)$, which is the subnetwork on a rate-$1$ Poisson point process $\Xi(1)$.  
Now augment the network $\SS(1)$ by including the point at the origin and the routes from the origin 
to each $\xi \in \Xi(1)$.  
The newly added edges have some random total  length $L$.
\begin{Proposition}
\label{Pmarginal}
$\Ex L = \ell/2$.
\end{Proposition}

\proof
Recall (\ref{ell-scale}) the scaling relation
 $\ell (\lambda) = \lambda^{1/2} \ell$, where 
  $\ell(\lambda)$ is the edge-intensity of the subnetwork $ \SS(\lambda)$ of a Poisson process of 
  point-intensity $\lambda$ . 
  Differentiating with respect to $\lambda$,
  \[ \ell^\prime(1) = \sfrac{1}{2} \ell . \] 
  So we need to show $\Ex L = \ell^\prime(1)$.

 Consider the space-time 
Poisson point process $(\Xi(\lambda), 0<\lambda <\infty)$ from section \ref{sec-TS}.
Each arriving point creates some additional network length, say 
$\tilde{L}(\xi)$, and for a point arriving at time $\lambda$,  write $\tilde{\ell}(\lambda)$ for the mean 
additional network length.  
Now
\[ \ell(\lambda_0) = \Ex \sum_{\xi \in \Xi(\lambda_0) \cap [0,1]^2} \tilde{L} (\xi) 
= \int_0^{\lambda_0} \tilde{\ell}(\lambda) \ d \lambda 
\]
and so $\ell^\prime (1) = \tilde{\ell}(1)$.

\subsection{A lower bound on network length}
Write $\Delta$ for the parameter $\Ex D_1$ of a SIRSN.
Write $\ell_*(\Delta)$ for the minimum possible value of $\ell$ in a SIRSN with a given value of $\Delta$.  

\begin{Proposition}
\label{PlD}
$\ell_*(\Delta) = \Omega((\Delta -1)^{-1/2})$ as $\Delta \downarrow 1$.
\end{Proposition}

The proof is based on a bound (Proposition \ref{PaL})  involving the geometry of deterministic paths, 
somewhat similar to bounds used in \cite{me116} section 4.
Figure 7 illustrates the argument to be used.

\setlength{\unitlength}{0.35in}
\begin{picture}(11,5)(-1,-1)
\put(-1,0){\line(1,0){11}}
\put(0,0){\line(0,1){0.2}}
\put(3,0){\line(0,1){0.2}}
\put(5,0){\line(0,1){0.2}}
\put(6,0){\line(0,1){0.2}}
\put(9,0){\line(0,1){0.2}}
\put(-0.2,-0.38){0}
\put(2.88,-0.38){L}
\put(4.84,-0.38){U}
\put(5.8,-0.38){2L}
\put(8.8,-0.38){3L}
\put(5,0){\line(0,1){3.3}}
\put(-1,0){\line(0,1){1}}
\put(0,0){\line(0,1){1}}
\put(-1,1){\line(1,0){1}}
\put(9,0){\line(0,1){1}}
\put(10,0){\line(0,1){1}}
\put(9,1){\line(1,0){1}}
\put(4,2){\line(2,1){3.2}}
\put(7.2,3.6){\line(2,-3){2.2}}
\put(9.4,0.3){\circle*{0.2}}
\put(4,2){\line(-1,1){1.2}}
\put(2.8,3.2){\line(-1,-2){1.8}}
\put(1,-0.4){\line(-3,1){1.8}}
\put(-0.8,0.2){\circle*{0.2}}
\put(-0.9,0.4){$z_1$}
\put(9.18,0.6){$z_2$}
\put(5,2.5){\line(1,0){1}}
\put(5.6,2.6){$\beta$}
\put(4.6,2.5){$\xi$}
\end{picture}

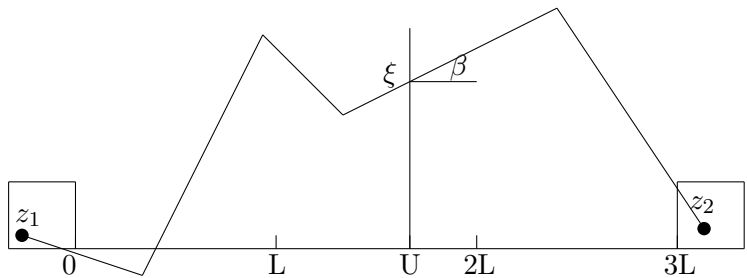
{\bf Figure 7.}

\begin{Proposition}
\label{PaL}
Let $\alpha, L, D$ and $\theta_0$ be positive reals satisfying 
$\theta_0< \pi/2$ and
\begin{eqnarray}
2 \sqrt{(D-1)^2 + (\sfrac{3}{2}L +1)^2} &=& (1+2\alpha) (3L+2) \label{D12}\\
L\left( \sfrac{1}{\cos \theta_0} -1 \right) &=& 4 \alpha \sqrt{(3L+2)^2 + 1} . \label{Lcos} 
\end{eqnarray}
Let $\RR$ be a route from 
some point $z_1$ in the unit square $[-1,0] \times [0,1]$ 
to some point $z_2$ in the unit square $[3L, 3L+1] \times [0,1]$, and suppose 
\begin{equation}
 \len(\RR) \le (1 + 2\alpha) |z_2 - z_1|  . 
 \label{Rlen}
 \end{equation}
Take $U$ uniform random on $[L,2L]$.  
The route $\RR$ first crosses the vertical line 
$\{(U,y), -\infty < y < \infty\}$ at some random point 
$(U,\xi(U))$ and at some angle 
$\beta(U) \in (- \sfrac{\pi}{2}, \sfrac{\pi}{2})$ relative to horizontal.   
Then \\
(i) $|\xi(U)| \le D $. \\
(ii) $\Pr (|\beta(U)| \le \theta_0) \ge \sfrac{1}{2} $.
\end{Proposition}
\proof
The maximum possible value of $\xi(U)$ arises in the case where
$z_1 = (-1,1), \ z_2 = (3L+1,1), \ U = \sfrac{3}{2}L$, the route consists of straight lines from 
$z_1$ to $(U,\xi(U))$ to $z_2$, and the route-length attains equality in (\ref{Rlen}).
In this case the value of $\xi(U)$ is the quantity $D$ satisfying (\ref{D12}), establishing (i).

Writing $\beta(u)$ for the angle (relative to horizontal) of the route at 
$x$-coordinate $u$,
then the length ($\Lambda$, say) of the route between $x$-coordinates $L$ and $2L$
equals 
$\int_L^{2L} \frac{1}{\cos \beta(u)}  \ du$.  
This implies
\[ \Lambda - L \ge (\sfrac{1}{\cos \theta_0 } - 1) \times L \Pr (\beta(U) \ge \theta_0)
. \]
But by considering excess length (relative to a horizontal route), (\ref{Rlen}) implies 
\[ \Lambda - L \le 2 \alpha |z_2 - z_1| \le 2 \alpha \sqrt{(3L+2)^2+1} . \]
Combining these inequalities gives a lower bound on $\Pr (\beta(U) \ge \theta_0)$ which equals $1/2$ 
when $\theta_0$ satisfies (\ref{Lcos}), establishing (ii).
\qed

\paragraph{Proof of Proposition \ref{PlD}}
Consider a SIRSN with parameters $\ell$ and $\Delta$ and with 
 induced subnetwork $\SS$ on a Poisson point process $\Xi$.
Set $\alpha = \Delta - 1$.
Suppose we can choose $L, D, \theta_0$ to satisfy, along with the given $\alpha$, the equalities  (\ref{D12},\ref{Lcos}) -- note this leaves us one degree of
freedom.

With probability $(1 - e^{-1})^2$ there are points $z_1$ and $z_2$ of the Poisson
process 
in the unit squares $[-1,0]\times [0,1]$ and $[3L,3L+1] \times [0,1]$.
By Markov's inequality and the definition of $\Delta$, with probability at least
$1/2$ the route 
$\RR(z_1,z_2)$ has length at most $(1+ 2\alpha)|z_2 - z_1|$.
Applying Proposition \ref{PaL} we deduce that, with probability 
$\ge (1-e^{-1})^2/4$, 
the network $\SS$ contains an edge that crosses the random vertical line 
$\{(U,y): \ - \infty < y < \infty\}$ at some point 
$(U,\xi(U))$ with 
$-D \le \xi(U) \le D$ 
and crosses at some angle $\beta(U) \in (-\theta_0, \theta_0)$ relative to horizontal.

If we translate vertically by $2D$, to consider routes between the unit squares 
$[-1,0]\times [2D,2D+1]$ and $[3L,3L+1] \times [2D,2D+1]$, then the potential
crossing points 
(using the same r.v. $U$) for the translated and untranslated cases are distinct.   
Now by considering translates by all multiples of $2D$, and noting that the
distribution of crossings 
of the random vertical line $\{(U,y): \ - \infty < y < \infty\}$ is the same as for
the 
$y$-axis, we have shown 
\begin{quote}
the mean intensity of crossings of the network $\SS$ over the $y$-axis at
angles 
$\in (-\theta_0, \theta_0)$ relative to horizontal is at least 
$\frac{(1-e^{-1})^2}{8D}$.  
\end{quote}
The stochastic geometry  identities (\ref{intersection-identity}, \ref{angle}) relates this mean intensity to the parameter
$\ell$ via 
\[ \mbox{this mean intensity} = \sfrac{\ell}{\pi} \int_{- \theta_0}^{\theta_0} \cos
\theta \ d \theta 
\le \sfrac{2 \ell \theta_0}{\pi} . \]
Combining with the previous inequality we find 
\[ \ell \ge \frac{1}{21 D \theta_0} . \]

Now set $L = \alpha^{-1/2}$ and consider the solutions of (\ref{D12},\ref{Lcos}) in the limit as $\alpha
\downarrow 0$: we find that solutions exist with
\[ \theta_0 \sim \sqrt{24 \alpha}; \quad  D \to 10 \]
which establishes Proposition \ref{PlD}.

\subsection{The minimum value of $\ell$ and the Steiner tree constant}
\label{sec-Steiner}
Take  $k$ uniform random points $Z_1,\ldots,Z_k$ in a square of area $k$ 
and consider the length $L_{ST}(k)$ of 
 the Steiner tree (the minimum-length connected network) on $Z_1,\ldots,Z_k$.
 Well-known subadditivity arguments \cite{steele97,yukich-book} 
imply that 
$\Ex L_{ST}(k) \sim \cst k$  for some constant $0< \cst < \infty$.  
One can define $\cst$ equivalently (see \cite{MR2235175} for results of this kind) as 
the infimum of $c$ such that there exists a translation-invariant connected
random network over 
$\Xi(1)$ with 
edge-intensity  $ c$. 
 From the latter description it is obvious that in any SIRSN we have $\ell \ge \cst$.
So the overall infimum
\begin{equation}
\ell_* := \mbox{ infimum of $\ell$ over all SIRSNs} 
\end{equation}
satisfies $\ell_* \ge \cst$, and below we outline an argument that 
the inequality is strict.
First we derive some simple lower bounds on $\cst$ and $\ell_*$.

(i) Write  $b(\xi)$ for the distance from $\xi$ to its closest neighbor in $\Xi(1)$.  
The discs of center $\xi$ and radius $b(\xi)/2$ are disjoint as $\xi$ varies and must contain
network length at least
$b(\xi)/2$, so 
\[ \cst \ge  \sfrac{1}{2} \Ex b(\xi) = \sfrac{1}{4}. \]   

(ii) In a network of edge-intensity $c$, (\ref{intersection-identity}) shows the mean number of edges 
crossing $\circl(0,r)$ equals 
$2 \pi r \times 2 \pi^{-1} c = 4rc$. 
If there is a point of $\Xi(1)$ inside $\disc(0,r)$ then there must be some such
crossing edge, so
\[ 1 - \exp(-\pi r^2) \le 4rc . \]
So
\[ \cst \ge \sup_r \frac{1 - \exp(-\pi r^2)}{4r} \approx 0.283 . \]

(iii) We can get a better bound on $\ell_*$ by using Proposition \ref{Pmarginal} as follows. 
Using the intensity calculation above, in a network of  edge-intensity $\ell$ the
probability that no 
edge crosses $\circl(0,r)$ is at least $1 - 4r \ell$.  
When a new point arives at $\xi$ in the $\SS(\lambda)$ process at time $\lambda
= 1$, if 
no existing edges cross $\circl(\xi,r)$ then the added network length $L$ is at
least $r$.
So 
\[ \Ex L \ge \sup_r r (1 - 4r \ell) = \sfrac{1}{16 \ell} . \]
But Proposition \ref{Pmarginal} says $\ell = 2 \Ex L$ and so we have shown
\begin{equation}
 \ell_* \ge \sqrt{1/8} \approx 0.353 . 
 \label{l*lb}
 \end{equation}
One could no doubt obtain small improvements by similar arguments. 

Here is an outline argument that $\ell_* > \cst$. \\
(i) In the Steiner tree on the Posson point process $\Xi(1)$, vertices of degree $>1$ have 
non-zero density, and their edges meet at some varying angles, whereas at the Steiner points (non-vertex junctions)  edges must meet at 120 degree angles. \\
(ii) If there were a SIRSN with $\ell \approx \cst$, then $\SS(1)$ would
 have essentially the properties (i).  But then in $\SS(1/2)$, obtained by deleting half the vertices of $\Xi(1)$ to get $\Xi(1/2)$, some of the deleted vertices would remain as junction points.   The ``varying angles" property 
implies the edge-intensity $\ell(1/2)$ of $\SS(1/2)$ is strictly larger than that of the Steiner tree on $\Xi(1/2)$, contradicting the scale-invariance property 
that the edge-intensities of $\SS(1)$ and the Steiner tree on $\Xi(1)$ 
are essentially equal.

\section{ General SIRSNs and their properties}
\label{sec-stronger}
In this section we study properties of $\SS(\lambda)$ in the $\lambda \to \infty$ limit, for a general SIRSN.  
Roughly speaking, this is studying ``the whole SIRSN" instead of sampled subnetworks, and such results 
depend on assumption (\ref{p1intro}).

Recall again the setup from (\ref{EE-def}) - (\ref{p1intro}).  So
$p(\lambda,r)$ is the edge-intensity of $\EE(\lambda,r)$, which is
the process of  points $z$ in edges of $\SS(\lambda)$ such that 
$z$ is in the route $\RR(\xi,\xi^\prime)$ for some 
 $\xi, \xi^\prime$ in $\Xi(\lambda)$ such that $\min(|z-\xi|, |z-\xi^\prime|) \ge r$.
Recall also from (\ref{p-scale})  the scaling relation 
$p(\lambda,r) = \lambda^{1/2} p(1,r\lambda^{1/2})$. 
Defining
\begin{equation}
p( r) := \lim_{\lambda \to \infty} p(\lambda, r) < \infty 
\label{strong-1}
\end{equation}
the assumption (\ref{p1intro}) that $p(1) < \infty$ and 
scaling  imply 
\begin{equation}
 p(r) = p(1) \times r^{-1}, \quad 0 < r < \infty . 
 \label{pr-scale}
 \end{equation}

\subsection{A connectivity bound}
Assumption (\ref{p1intro}) has a direct implication for the qualitative structure of a SIRSN: 
all the routes linking two regions, once they get away from a neighborhood of the regions, use only a finite number of different paths.
We first give a version of this result in terms of discs.  

\setlength{\unitlength}{0.1in}
\begin{picture}(20,19)(-20,-9)
\put(0,0){\circle{5}}
\put(0,0){\circle{10}}
\put(0,0){\circle{15}}
\put(1.3,0.9){\line(2,1){7}}
\put(0.4,1.6){\line(1,0){2.22}}
\put(2.1,0.2){\line(2,3){1.1}}
\put(-4,-4){\line(1,1){3}}
\put(-4,-4){\line(-1,0){3.6}}
\put(-5,-4){\line(-1,-3){0.8}}
\put(-5.2,-4.6){\line(-3,-1){1.3}}
\put(-2.2,-2.2){\line(4,1){1.7}}
\put(-2.4,-2.4){\line(1,3){0.8}}
\put(5.9,3.2){\line(1,0){2}}
\put(5.7,3.1){\line(1,4){0.7}}
\end{picture}

{\bf Figure 8.}
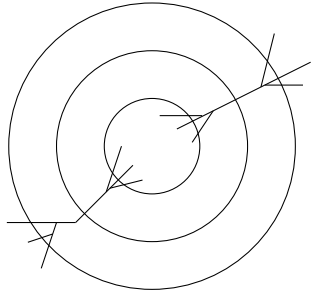
{\small Schematic for routes from inside $\disc(0,1/2)$ to outside $\disc(0,3/2)$ crossing the unit circle.}

\begin{Proposition}
\label{P-strong-1}
Take $0<r<1$ and
let $N(\lambda, r)$ be the number of distinct points on the unit circle at which 
some route $R(\xi,\xi^\prime)$ between some $\xi \in \Xi(\lambda) \cap \disc(0,1-r)$ 
and some $\xi^\prime \in \Xi(\lambda) \cap (\Reals^2 \setminus \disc(0,1+r))$ crosses the unit circle.  
Then
\[ \Ex \lim_{\lambda \to \infty} N(\lambda,r) \le 4 p(1) \ r^{-1} . \]
\end{Proposition}

\proof
Any crossing point is in $\EE(\lambda,r)$ and so by identity (\ref{intersection-identity})
\[ \Ex N(\lambda,r) \le 2\pi \times 2 \pi^{-1} p(\lambda,r) < \infty \]
and the result follows from (\ref{pr-scale}).
\qed

The following general version can be proved similarly.
\begin{Proposition}
 Let $\eps > 0$ and 
let $K_1, K_2$ be compact sets whose $\eps$-neighborhoods 
$K_1^\eps, K_2^\eps$ are disjoint.  
For $z_1 \in K_1, z_2 \in K_2$ let 
$\RR_\eps(z_1,z_2)$ be the subroute of $\RR(z_1,z_2)$ 
crossing from the boundary of $K_1^\eps$ to the boundary of 
 $K_2^\eps$. 
Let $N(\lambda, \eps, K_1, K_2)$ be the number of distinct paths amongst  the set
$\{\RR_\eps(\xi_1,\xi_2)\  : \ \xi_i \in \Xi(\lambda) \cap K_i\}$. Then 
\[\Ex \lim_{\lambda \to \infty} N(\lambda, \eps, K_1, K_2)< \infty . \]
\end{Proposition}

\subsection{A bound on normalized length}

\begin{Lemma}
\label{Lellp}
$\ell \leq 2 p(1)$.
\end{Lemma}
\proof
Define
\[ \RR_{\delta}(\xi,\xi^\prime) = 
\RR(\xi,\xi^\prime) \cap (\disc(\xi,\delta) \cup \disc(\xi^\prime,\delta)) \]
in words, the part of the route that is {\em within} distance  $\delta$  from one or both endpoints.
Then define
\[ \widehat{\SS}_{\delta}(\lambda) = \cup_{\xi, \xi^\prime \in \Xi(\lambda)} 
\RR_{\delta}(\xi,\xi^\prime) . \]
Note that clearly
\begin{equation}
S(\lambda) \setminus \widehat{\SS}_1(\lambda) \subseteq \EE(\lambda,1) .
\label{ShatS}
\end{equation}
By considering $\lambda = 1$,
\[ \ell \le p(1,1) + \iota (\widehat{\SS}_1(1) ) \]
where $\iota(\cdot)$ denotes edge-intensity.
Now write 
\[ \iota (\widehat{\SS}_1(1) ) = \sum_{k \ge 1} 
 \iota (\widehat{\SS}_{2^{1-k}}(1) \setminus \widehat{\SS}_{2^{-k}}(1)  ) . \]
For fixed $k \ge 1$, scaling by $2^k$ gives 
\begin{eqnarray*}
 \iota (\widehat{\SS}_{2^{1-k}}(1) \setminus \widehat{\SS}_{2^{-k}}(1)  ) 
&=& 2^{-k} \iota (\widehat{\SS}_{2}(2^{-2k}) \setminus \widehat{\SS}_{1}(2^{-2k})  ) \\
&\le & 2^{-k} p(2^{-2k}, 1) \mbox{ by (\ref{ShatS})} .
\end{eqnarray*}
So 
\[ \ell \le \sum_{k \ge 0} 2^{-k} p(2^{-2k}, 1) \le \sum_{k \ge 0} 2^{-k} p( 1)
. \]
\qed

\subsection{The network $\EE(\infty,r)$ of major roads}
Intuitively, the point of assumption (\ref{p1intro}) and the scaling 
relation (\ref{pr-scale}) is that we can define a proces 
$\EE(\infty, r):= \cup_{\lambda < \infty} \EE(\lambda,r)$ 
which must have edge-intensity 
$p(r) = p(1)/r$, and that in results like Proposition \ref{P-strong-1} 
we can replace $\lim_{\lambda \to \infty} N(\lambda, r)$ by $N(\infty,r)$.  
We don't want to give details of a completely rigorous treatment, but let us just suppose 
we can set up 
 $\EE(\infty,r)$ as a 
random element of some suitable measurable space, 
as we did for $\SS(\lambda)$ in section \ref{sec-TS-2}.

The conceptual point is that $\SS(\lambda)$ and $\EE(\lambda,r)$ depend on the external randomization, that is on the fact that we 
were studying a SIRSN via the random points
$\Xi(\lambda)$, 
but as outlined below $\EE(\infty,r)$ doesn't depend on such  external randomization.
Intuitively this is simply because $\cup_\lambda \Xi(\lambda)$ is dense in $\Reals^2$; 
we outline a measure-theoretic argument below.
\begin{Proposition}
\label{P*}
The FDDs $(\span(z_1,\ldots,z_k))$ of a SIRSN 
can be extended to a joint distribution, of these FDDs jointly with a random process $\EE^*(\infty,r)$, such
that, 
for any space-time PPP $(\Xi(\lambda), 0<\lambda<\infty)$ independent of the FDDs,
we have 
$\EE^*(\infty, r):= \cup_{\lambda < \infty} \EE(\lambda,r)$ a.s.
\end{Proposition}
{\bf Outline proof.}
For a suitable formalization of ``random subset of $\Reals^2$" we have the implication
\begin{quote}
if $\AA_1$ and $\AA_2$ are i.i.d. random subsets, and if 
$\AA_1 \cup \AA_2 \subseteq_{a.s.} \AA^\prime \ed \AA_1$, then 
$\AA_1 = A$ a.s. for some non-random subset $A$
\end{quote}
and then the corresponding ``conditional" implication
\begin{quote}
if $Z$ is a random element of some space, if $\AA_1$ and $\AA_2$ are random subsets 
conditionally i.i.d. given $Z$, and if 
$\AA_1 \cup \AA_2 \subseteq_{a.s.} \AA^\prime$ where 
$(Z,\AA^\prime) \ed (Z,\AA_1)$, then 
$\AA_1 = \AA$ a.s. for some $Z$-measurable random subset $A$.
\end{quote}
So take two independent space-time PPPs 
$\Xi^1(\lambda), \Xi^2(\lambda)$ 
and use a measure-preserving bijection 
$[0, \infty) \cup [0,\infty) \to [0,\infty)$ to define another space-time PPP $\Xi^\prime(\lambda)$ 
in terms of $\Xi^1$ and $\Xi^2$.  
The associated networks satisfy 
\[ \EE^1(\infty,r) \cup \EE^2(\infty,r) = \EE^\prime(\infty,r) \ed \EE^1(\infty,r) \] 
and this holds jointly with the FDDs of the SIRSN.  
Since $ \EE^1(\infty,r)$ and $\EE^2(\infty,r)$ are conditionally i.i.d. given the SIRSN. 
Proposition \ref{P*} follows from the general ``conditional implication" above.

\subsection{Transit nodes and shortest path algorithms}
\label{sec-transit}
Here we make a connection with the ``shortest path algorithms" literature 
mentioned in section \ref{sec-intro-SPA}.

Fix $h$ and take the square grid of lines with inter-line spacing equal to $h$.  
Define $\TT_h$ to be the set of points of intersection of 
$\EE(\infty,h)$ with that grid. 
\begin{Lemma}
\label{LTh}
(i) $\TT_h$  has point-intensity   $4 \pi^{-1} h^{-2} p(1)$. \\
(ii) For each $z \in \Reals^2$ there is a subset $T_z$ of $\TT_h$, of mean size 
$24 \pi^{-1} p(1)$, and with $|z^\prime - z| \le 2^{3/2} h$ for each $z^\prime \in T_z$, 
such that for each pair $z_1, z_2$ with $|z_2 - z_1| > 3h$ 
the route $\RR(z_1,z_2)$ passes through some point of $T_{z_1}$ and some point of $T_{z_2}$.
\end{Lemma}
\proof 
The grid has edge-intensity $2h^{-1}$, so from
(\ref{intersection-identity}) the point-intensity of $\TT_h$ is
$2\pi^{-1} \times p( h) \times 2h^{-1}$, 
and (i) follows from scaling  (\ref{pr-scale}).

For any starting point $z$ consider the closest grid intersection 
$(ih,jh)$.
Then $z$ is in some square with corner $(ih,jh)$, 
say the square $[(i-1)h,ih] \times [jh,(j+1)h]$.   
Let $T_z$ be the set of points of intersection of 
$\EE(\infty,h)$ with the concentric square 
$S_z = [(i-2)h, (i+1)h] \times [(j-1)h, (j+2)h]$. 
This square has boundary length $12h$ and so the mean size of $T_z$ 
equals $2\pi^{-1} \times p( h) \times 12h = 24 \pi^{-1} p(1)$.  
By construction 
\[ \sfrac{3}{2} h <  |z^\prime - z| \le 2^{3/2} h \mbox{  for each } z^\prime \mbox{ on the boundary of}  S_z  \]
and in particular for each $z^\prime \in T_z$. 
If $|z_2 - z_1| > 3h$ then the squares $S_{z_1}$ and $S_{z_2}$ do not overlap, and the points 
$z^\prime_1$ and $z^\prime_2$ at which the route crosses their boundaries are in $\TT_h$.
\qed

\paragraph{Informal algorithmic implications}
One cannot rigorously relate our ``continuum" setup to discrete algorithms, 
but in talks we present the following informal calculation.
For the real-world road network in a country we have empirical statistics
\begin{itemize}
\item $A$: area of country
\item $\eta$: average number of road segments per unit area
\item $M  = \eta  A$: total number of road segments in country
\item $p(r)$: ``length per unit area" of the subnetwork consisting of 
 segments on routes with start/destination each at distance $>r$ from the segment. 
\end{itemize}
For a real-world network   
there is an inconsistency between scale-invariance and having a finite number $\eta$ 
of road segments per unit area,
but let us imagine approximate scale-invariance over scales of say 2 - 100 miles, 
and modify a scale-invariant model by deleting road segments of very short length.  
In what follows it is helpful to imagine the unit of length to be (say) 20 miles.
 
Fix $r$.  Lemma  \ref{LTh} (with $h = r$) suggests that in the real-world network we can find transit nodes such that there are 
 $O(p(1))$ transit nodes within distance $O(r)$ of a typical point.  
If so then we can analyze
the  algorithmic procedure outlined in section  \ref{sec-intro-SPA}.
The local search involves a region of radius $r$ and hence with $O(\eta r^2)$ edges.
Regarding the time-cost of a single Dijkstra search as $c_1 \times (\mbox{number of edges})$, 
the time-cost of finding the route to each local transit node is 
$O\left(c_1  \ (\eta r^2) p(1) \right)$.  
Transit nodes have point-intensity $O(p(1)/r^2)$, so the total number is $O(Ap(1)/r^2)$. 
Regard the space-cost of storing a $k \times k$ matrix of inter-transit-node routes as $c_2 k^2$; so this space-cost is $O \left( c_2 \ (p(1)A/r^2)^2 \right)$. 
Summing  the two costs and optimizing over $r$,  the
optimal cost is 
$O(c_1^{2/3}c_2^{1/3}\eta^{2/3}A^{2/3}p^{4/3}(1))
= O(c_1^{2/3}c_2^{1/3}p^{4/3}(1)M^{2/3})$ 
and this $O(M^{2/3})$ scaling represents the improvement over the $O(M)$  scaling for Dijkstra.  
The corresponding optimal number of transit nodes is 
$O((c_1/c_2)^{1/3} p^{2/3}(1) M^{1/3})$.  
The latter has a more interpretable formulation. 
If the only alternative algorithms were a Dijkstra search of cost  $c_1 \times (\mbox{number of edges})$ 
or table look-up of cost $c_2 \times (\mbox{number of edges})^2$ , then there would be some critical  number of
edges  at which 
one should switch between them, and this is just the solution $\mcrit$ of $c_1 \mcrit =
c_2 \mcrit^2$.  So the optimal number of transit nodes is  $O(\mcrit^{1/3}
p^{2/3}(1) M^{1/3})$.

\subsection{Number of singly infinite geodesics}

Write $\SS^*(\lambda)$  for the spanning subnetwork obtained from $\SS(\lambda)$  by adding a city at the origin $\origin$.  This process inherits the scaling-invariance property (\ref{SS-def-4}) of $\SS(\lambda)$.
We know from (\ref{Lsingly}) that at least one singly infinite geodesic from $\origin$ exists. 
The set of all singly infinite geodesics in $\SS^*(\lambda)$ from $\origin$ forms {\em a priori} a tree, 
because two geodesics that branch cannot re-join, by route compatability property (iv) from section \ref{sec-TS}.
So consider 
\[ q(\lambda,r) := \Ex 
\mbox{(number of distinct points at which some singly infinite geodesic } \]
\vspace*{-0.3in}
\[ \mbox{  in $\SS^*(\lambda)$ from $\origin$ first crosses the circle of radius $r$).}
\]
What we know in general is
\[ 1 \le q(\lambda, r) \le \infty; \quad
r \to q(\lambda,r) \mbox{ is increasing;} \quad 
\lambda \to q(\lambda,r) \mbox{ is increasing}  \] 
and the scaling property gives
\begin{equation}
 q(\lambda,r) = q(1,r\lambda^{1/2}) . 
 \label{q-scale}
 \end{equation}
So the $\lambda \to \infty$ limit 
$q(\infty,r) := \lim_{\lambda \to \infty} q(\lambda,r)$ exists (maybe infinite), and the scaling property implies
\[ q(\infty,r) = q(\infty,1) \in [1,\infty], \quad 0 < r < \infty . \]
So consider the property 
\begin{equation}
q(\infty,1) < \infty .
\label{strong-2}
\end{equation}
By applying Proposition \ref{P-strong-1} with $r \approx 1$ we see 
\begin{equation}
q(\infty,1) \le 4 p(1) .
\end{equation}
So (\ref{strong-1}) implies (\ref{strong-2}).
So we have shown the following.
\begin{Corollary}
\label{CG}
As $\lambda \to \infty$ the number  of singly infinite geodesics in $\SS^*(\lambda)$ from $\origin$ 
increases to a finite limit number (perhaps a random number with finite mean) $G$. 
Moreover, if $G>1$ then these geodesics branch at $\origin$.
\end{Corollary}

\section{Unique singly-infinite geodesics and continuity}
\label{sec-cul}
For a SIRSN, let us call the property $G = 1$ a.s. (in the notation of 
Corollary \ref{CG} above) the 
{\em unique singly-infinite geodesics} property.  
It is conceivable that this property always holds -- we record this later in
Open Problem \ref{OP-counter}.  
Uniqueness of geodesics is closely related to continuity of routes $\RR(z_1,z_2)$ as $(z_1,z_2)$ vary, as will be seen 
in section \ref{sec-continuity}.

\subsection{Equivalent properties}
Here we show that several properties, the simplest being (\ref{cul-1}), are equivalent to the 
unique singly-infinite geodesics property. We will give definitions and proofs as we proceed, and then summarize as Proposition \ref{Pcul}.

Consider two independent uniform random points $U_1, U_2$ in $\disc(\origin, 1)$.  
By the route-compatability property, the intersection of 
$\RR(\origin, U_1)$ and $\RR(\origin, U_2)$ is a sub-route from 
$\origin$ to some {\em branchpoint} $B_{1,2}$, where either $B_{1,2} \neq \origin$ 
or the intersection consists of the single point $\origin$ 
(in which case, set $B_{1,2} = \origin$).  
So we can define a property
\begin{equation}
\Pr (B_{1,2} = \origin) = 0 .
\label{cul-1}
\end{equation}

\paragraph{Unique singly-infinite geodesics imply (\ref{cul-1}).} 
Suppose (\ref{cul-1}) fails. 
Then there exists $\eps > 0$ such that, for independent random points 
$U^1_1, U^1_2$ in $\disc(\origin, 1) \setminus \disc(\origin, \eps) $,
their branchpoint $B^1_{1,2}$ satisfies 
$\Pr (B^1_{1,2} = \origin) \ge \eps$. 
Scaling by $\eps^{-m}, m \ge 1$ and using scale-invariance, 
for independent random points 
$U^m_1, U^m_2$ in $\disc(\origin, \eps^{-m}) \setminus \disc(\origin, \eps^{1-m}) $,
their branchpoint $B^m_{1,2}$ satisfies 
$\Pr (B^m_{1,2} = \origin) \ge \eps$.  
It follows that, with probability $\ge \eps - o(1)$ as $m \to \infty$, 
there exists points $\xi^m_1, \xi^m_2$ of 
$\Xi(1) \cap ( 
\disc(\origin, \eps^{-m}) \setminus \disc(\origin, \eps^{1-m}) 
)$ 
such that
\[ \mbox{ routes $\RR(\origin, \xi^m_1)$ and $\RR(\origin, \xi^m_2)$ 
branch at $\origin$}. \]
So on an event of probability $\ge \eps$ this property holds for infinitely many $m$.  Then on that event we have $G > 1$,  by compactness within the spanning subnetwork 
$\SS^*(1)$ (Lemma \ref{Lss1}). 
\qed

\medskip \noindent
Next consider the spanning subnetwork 
$\SS^*(\lambda)$ on points $\Xi(\lambda) \cup \{\origin\}$.  
The intersection of all routes 
$\RR(\origin, \xi), \ \xi \in \Xi(\lambda)$ 
is a sub-route from 
$\origin$ to some branchpoint $B(\lambda)$.
So we can define a property
\begin{equation}
\Pr(B(1) = \origin) = 0 .
\label{cul-2}
\end{equation} 
Clearly (\ref{cul-2}) implies (\ref{cul-1}); we need to argue the converse.

\paragraph{(\ref{cul-1}) implies (\ref{cul-2}).} 
Suppose (\ref{cul-1}).
For each $r< \infty$ the intersection of routes 
$\RR(\origin, \xi), \ \xi \in \Xi(1) \cap \disc(\origin, r)$ is a subroute  
$\pi(1,r)$ from $\origin$ to some branchpoint $B(1,r)$, 
and by (\ref{cul-1}), scaling and the finiteness of $\Xi(1) \cap \disc(\origin, r)$ 
we have 
\begin{equation}
\Pr ( B(1,r) = \origin) = 0 , \mbox{ each } r < \infty .
\label{Bir}
\end{equation}
As $r$ increases the subroute $\pi(1,r)$ can only shrink, and the quantity in 
(\ref{cul-2}) is the limit 
$B(1) = \lim_{r \to \infty} B(1,r)$. 
To prove (\ref{cul-2}) it suffices, by (\ref{Bir}), to prove 
\begin{equation}
 B(1,r) \mbox{ is constant for all large $r$, a.s. }
\label{B1r}
\end{equation}
We may suppose (otherwise the result is obvious) that for some $r_0 \ge 4$ 
the subroute $\pi(1,r_0)$ stays within $\disc(\origin, 1)$.  
As $r$ increases, the only way that $B(1,r)$ can change at $r$ is if there is a 
point $\xi \in \Xi(1) \cap \circl(\origin, r)$ for which the route 
$\RR(\origin,\xi)$ diverges from the existing subroute $\pi(1,r-)$ before the existing 
branchpoint $B(1,r-)$.  
If this happens, at $r_1$ say, then consider the subroute 
$\theta(r_1) = \RR(\origin,\xi) \cap(\disc(\origin, 4) \setminus \disc(\origin, 1))$ 
which has length at least $3$. 
Now suppose $B(1,r)$ again changes at some larger value $r_2$. 
Then the corresponding subroute $\theta(r_2)$ must be disjoint from $\theta(r_1)$, 
by the route-compatability property. 
Now the ``finite length in bounded regions" property (\ref{BLBA}) implies that 
$B(1,r)$ can change at only finitely many large values of $r$, 
establishing (\ref{B1r}).
\qed

\medskip \noindent
Now make a slight re-definition of $B(\lambda)$, by considering only points $\xi$ outside the unit disc.
That is, 
the intersection of all routes 
$\RR(\origin, \xi), \ \xi \in \Xi(\lambda)\setminus \disc(\origin, 1) $ 
is a sub-route $\tilde{\pi}(\lambda)$ from 
$\origin$ to some branchpoint $B_1(\lambda)$. 
Using scale-invariance it is easy to check that (\ref{cul-2}) is equivalent to 
\begin{equation}
 \Pr (B_1(\lambda) = \origin) = 0 \mbox{ for each } \lambda < \infty . 
 \label{B1L}
 \end{equation}
As $\lambda$ increases, the sub-routes $\tilde{\pi}(\lambda)$
can only shrink, and the intersection of these subroutes over all $\lambda < \infty$  
is again a subroute from $\origin$ to some point $B_1(\infty)$.  
So we can define a property
\begin{equation}
\Pr (B_1(\infty) = \origin) = 0 .
\label{cul-3}
\end{equation} 
Clearly (\ref{cul-3}) implies (\ref{B1L}); we need to argue the converse.

\paragraph{(\ref{B1L}) implies (\ref{cul-3}).} 
Suppose (\ref{B1L}).  
To prove (\ref{cul-3}) we essentially repeat the argument above, but use assumption 
(\ref{p1intro}) instead of (\ref{BLBA}).  
It is enough to show that, as $\lambda$ increases,  $B_1(\lambda)$ can change at only finitely many large 
values of $\lambda$.  
And we may suppose that for large $\lambda$ 
the subroute $\tilde{\pi}(\lambda)$ stays within $\disc(\origin, 1/4)$.  
If $B_1(\lambda)$ changes at $\lambda_1$ then there is a point $\xi$ appearing at 
``time" $\lambda_1$ for which $\RR(\origin, \xi)$ diverges from the existing subroute 
$\tilde{\pi}(\lambda_1 -)$ and so must cross $\circl(\origin, 5/8)$ at some point 
$z(\lambda_1) \in \EE(\lambda_1, 3/8) \subset \EE(\infty, 3/8)$.  
By route-compatability the points $z(\lambda_i)$ corresponding to different values 
$\lambda_i$ where $B_1(\lambda)$ changes must be distinct, and then (\ref{p1intro}) 
implies $\EE(\infty, 3/8) \cap \circl(\origin, 5/8)$ is an a.s. finite set of points.
\qed

\medskip \noindent
Clearly (\ref{cul-3}) implies unique singly-infinite geodesics, by the final assertion of Corollary \ref{CG}.  
We have now shown a cycle of equivalences. 
Finally, by scaling (\ref{cul-3}) is equivalent to the following property, where the notation is chosen to be consistent with 
notation in the next section.  
Define $Q(\lambda, 0, B)$ to be the probability that the routes 
$\RR(\origin,\xi^\prime)$ to all  points $\xi^\prime \in \SS(\lambda) \cap (\disc(\origin,B))^c $ 
do {\bf not} all first exit $\disc(\origin, 1)$ at the same point.  
Then (\ref{cul-3}) is equivalent to 
\begin{equation}
\lim_{B \uparrow \infty} \lim_{\lambda \to \infty} Q(\lambda, 0, B) = 0 .
\label{def-weak-uniq}
\end{equation}
To summarize:
\begin{Proposition}
\label{Pcul}
Properties (\ref{cul-1}), (\ref{cul-2}), (\ref{Bir}), (\ref{cul-3}) and
(\ref{def-weak-uniq}) are each equivalent to the 
unique singly-infinite geodesics property.
\end{Proposition}

\subsection{Continuity properties}
\label{sec-continuity}
In the previous section we studied properties of long routes from a single point. 
We now consider long routes from nearby points, and in this context it seems harder to understand whether different properties are equivalent. 
Suppose, for this discussion, the unique singly-infinite geodesics property 
holds.  Then the geodesics from $\origin$ and from 
$\one = (1,0) \in \Reals^2$ are either disjoint or coalesce; we do not know (Open Problem \ref{OP-counter}) whether the property 
\begin{equation}
\mbox{ 
 the geodesics from $\origin$ and from 
$\one$ coalesce a.s }
\label{prope01}
\end{equation}
always holds or is stronger.  There are several equivalent ways of saying 
(\ref{prope01}) -- see the end of this section -- but what's relevant now is that it is equivalent to the property that, for each $\lambda$, the geodesics from each point of 
$\SS^*(\lambda) \cap \disc(\origin,1)$ coincide outside a disc of random radius $R(\lambda) < \infty$ a.s..  So we can then ask whether the 
property
\[ R(\infty) := \lim_{\lambda \to \infty} R(\lambda) < \infty \mbox{ a.s. } \]
is implied by property (\ref{prope01}) or is stronger.  We restate this latter property as 
(\ref{def-strong-uniq}) below.

For $0 < \eps < 1 < B$ define 
$Q(\lambda, \eps, B)$ to be the probability that the routes 
$\RR(\xi,\xi^\prime)$ between points $\xi \in \SS(\lambda) \cap \disc(\origin,\eps)$
and points $\xi^\prime \in \SS(\lambda) \cap (\disc(\origin,B))^c $ 
do {\bf not} all first exit $\disc(\origin, 1)$ at the same point.  
Note $Q(\lambda, \eps, B)$ is monotone increasing at $\lambda$ increases, 
and decreasing as $B$ increases or $\eps$ decreases. 
So we can define 
\[ Q(\infty, \eps, B) := \lim_{\lambda \to \infty}   Q(\lambda, \eps, B) \]
and then define a property of a SIRSN
\begin{equation}
\lim_{\eps \downarrow 0, B \uparrow \infty} Q(\infty, \eps, B) = 0 
\label{def-strong-uniq}
\end{equation}
where the limit value is unaffected by the order of the double limit. 
In words, (\ref{def-strong-uniq}) says that (with high probability) every route from a small neighborhood of the origin 
to any distant point will first cross the unit circle at the same place.   
Property (\ref{def-strong-uniq}) implies (\ref{prope01}) and implies form (\ref{def-weak-uniq}) 
of the unique singly-infinite geodesics property, 
which is the same assertion for routes from the origin only.

The kinds of properties described above relate to questions about
continuity of the routes 
$\RR(z_1,z_2)$ as $z_1, z_2$ vary, and we will give one such relation as  
Lemma \ref{L-strong-cont} below.

Consider $0 < \eta < \delta < 1/2$ and for points  
$\xi \in \SS(\lambda) \cap \disc(\origin,\eta)$ and 
$\xi^\prime \in \SS(\lambda) \cap \disc(\one,\eta)$ 
with route $\RR(\xi,\xi^\prime)$ let $\RR_\delta(\xi,\xi^\prime)$ 
be the sub-route between the first exit from $\disc(\origin,\delta)$ and 
the last entrance into $\disc(\one,\delta)$.  
Let $\Psi(\lambda, \eta, \delta)$ be the probability that the sub-routes 
$\RR_\delta(\xi,\xi^\prime)$ for all 
$\xi \in \SS(\lambda) \cap \disc(\origin,\eta)$ 
and all 
$\xi^\prime \in \SS(\lambda) \cap \disc(\one,\eta)$ are {\bf not} 
all an identical sub-route.  
As above, by monotonicity we can define 
\[ \Psi(\infty, \eta, \delta) := \lim_{\lambda \to \infty}   \Psi(\lambda, \eta, \delta) \]
and then define a property of a SIRSN
\begin{equation}
\lim_{\eta \downarrow 0} \Psi(\infty, \eta, \delta) = 0 \quad \forall \delta .
\label{def-strong-cont}
\end{equation}
In words, (\ref{def-strong-cont}) says that (with high probability) all routes from a very small neighborhood of the origin 
to a very small neighborhood of $\one$ coincide outside of larger small neighborhoods.

\begin{Lemma}
\label{L-strong-cont}
Property (\ref{def-strong-uniq}) implies 
property (\ref{def-strong-cont}).
\end{Lemma}
\proof 
Choose $a$ such that $a \eta < 1 < a \delta$.
Take the definition of $\Psi$, scale by $a$, and use scale-invariance 
to obtain the following.
\begin{quote}
The probability that the sub-routes 
$\RR_{a \delta}(\xi,\xi^\prime)$ for all 
$\xi \in \SS(a^{-2} \lambda) \cap \disc(\origin, a \eta)$ 
and all 
$\xi^\prime \in \SS(a^{-2} \lambda) \cap \disc((a,0), a \eta)$ are {\bf not} 
all an identical sub-route 
equals $\Psi(\lambda, \eta, \delta)$.
\end{quote} 
When this occurs there are two non-identical sub-routes between 
$\circl(\origin, a \delta)$ and $\circl( (a,0), a \delta)$, 
which imply two non-identical sub-routes between 
$\circl(\origin, 1)$ and $\circl( (a,0), 1)$. 
For this to happen, either the defining event for 
$Q(a^{-2} \lambda, a \eta, a/2)$, 
or the analogous event with reference to $(a,0)$ instead of $\origin$,
must occur; otherwise all routes in question pass through the same points 
on $\circl(\origin, 1)$ and $\circl( (a,0), 1)$, contradicting the  
route-compatability properties of section \ref{sec-TS}.  
So
\[ \Psi(\lambda, \eta, \delta) \leq 2 Q(a^{-2} \lambda, a \eta, a/2) . \]
Letting $\lambda \to \infty$
\[ \Psi(\infty, \eta, \delta) \leq 2 Q(\infty, a \eta, a/2) . \]
Choosing $a = \eta^{-1/2}$ 
establishes the lemma.
\qed

{\bf Remark.}  Lemma \ref{L-strong-cont} is almost enough to prove that, under condition (\ref{def-strong-uniq}),
we have the continuity property 
\begin{equation}
\mbox{ if } (z^n_1,z^n_2) \to (z_1,z_2) \mbox{ then } \RR(z^n_1,z^n_2) \to \RR(z_1,z_2) \mbox{ a.s. }
\label{eq-ct}
\end{equation}
where convergence of paths is in the sense of  section \ref{sec-TS-2}.  
To deduce (\ref{eq-ct}) one would need also to show that the lengths of 
$ \RR(z^n_1,z^n_2) \cap (\disc(z_1,\eps_n) \cup \disc(z_2,\eps_n))$ 
tend to $0$ a.s. for all $\eps_n \to 0$.   
This is loosely related to Open Problem \ref{OP-len}.

\paragraph{Another property equivalent to (\ref{prope01}).} 
Because geodesics either colalesce or are disjoint, for any countable set of initial points there is some set of  ``geodesic ends", where each such ``end" corresponds to a tree of coalescing geodesics from originating ``leaves".  
By a small modification of the proof of Corollary \ref{CG}, the mean number of 
such ends from the points $\Xi(\lambda) \cap \disc(\origin, 1)$ is at most 
$4 p(1)$, so we can let $\lambda \to \infty$ and deduce that the number 
$G^* \ge 1$ of ends from initial  points $\Xi(\infty) \cap \disc(\origin, 1)$  satisfies $\Ex G^* \le 4 p(1)$.  Then by scale-invariance, for each 
$0<r<\infty$ the number  of ends from initial  points $\Xi(\infty) \cap \disc(\origin, r)$ equals $G^*$.  So the property
\[ G^* = 1 \mbox{ a.s. } \]
is clearly equivalent to property (\ref{prope01}) (plus the 
unique singly-infinite geodesics property).  
Note that if $G^* > 1$ then there are
a finite number of different ``geodesic trees" each of whose leaf-sets is dense in $\Reals^2$ -- behavior hard to visualize.

\subsection{The binary hierarchy model}

\begin{Proposition}
The binary hierarchy model has property (\ref{def-strong-uniq}) .
\end{Proposition}
\proof 
Consider the last stages of construction of the model in 
section \ref{sec-1-complete}.  Rotation and scaling do not affect the property of interest, so it will suffice to prove the property
in the model $\RRti$.  Consider 
the argument from  ``proof of Proposition \ref{PnullA}" in section 
\ref{sec-1-plane} but with large rescalings
of $B$ instead of small rescalings.  Combining this argument with the  
construction of $\RRti$ at the start of section \ref{sec-1-complete} one can show (details omitted)  that the set 
\[ A^\prime:= \{z \in \Reals^2:  \ z \mbox{ in only finitely many } B^\prime_i\} \] has area zero; here  
$B^\prime_i:= \sigma_{2^{i}}B$ is the ``large" rescaling of 
the union $B: = \cup_G \Sigma_G$ of the translates $\Sigma_G$ of the 
small subsquare $\Sigma$ of the
basic $2^{h+1} \times 2^{h+1}$ square $G$ in the Figure 6 configuration.
By translation-invariance, this implies that a.s. $\origin \not\in A^\prime$.  
For such a realization there is a random infinite sequence $i(j)$ with 
$\origin \in \sigma_{2^{i(j)}}B$, and  any singly-infinite geodesic from $\origin$ must pass through 
the corresponding infinite sequence $b_{i(j)}$ of points determined by Figure 6.  This establishes the  unique singly-infinite geodesic property.
Moreover $\origin$ lies in some translated square $\Sigma_{i(j)}$ of side 
$\eps 2^{i(j)}$ and for any other point in that square its geodesic must coalesce with the geodesic from $\origin$ at or before  $b_{i(j)}$.  
It is easy to check that the squares $\Sigma_{i(j)}$ eventually cover any fixed disc, and this establishes property (\ref{def-strong-uniq}).
\qed

\section{Open problems and final discussion}
\label{sec-final-disc}

\subsection{Other specific models?}
\label{sec-OSM}
A major challenge is finding other explicit examples of SIRSN models.  
Let us pose the vague problems
\begin{OP}
\label{OP1}
Give a construction of a SIRSN which is ``mathematically natural" in some sense, e.g. in the sense that 
there is an explicit formula for the distribution of subnetworks $\span(z_1,\ldots,z_k)$.
\end{OP}
\begin{OP}
\label{OPx1}
Give a construction of a SIRSN which is ``visually realistic" in the sense of not looking very different from a real-world road network.
\end{OP}

\subsection{Quantitative bounds on statistics}
In designing a finite road network there is an obvious  tradeoff between total length and the network's effectiveness in providing short routes,  
so in our context there is a tradeoff between $\ell$ and $\Delta:= \Ex D_1$. 
More generally 
\begin{OP}
\label{OP-QBS}
What can we say about the set of possible values, over all SIRSNs, of the triple 
$(\Delta = \Ex D_1, \ell, p(1))$ of statistics of a SIRSN?  
\end{OP}
This is a sensible question because each statistic is dimensionless,
that is not dependent on choice of unit of length -- a non-dimensionless statistic 
would take all values in $(0,\infty)$ by scaling.  

We have given three results relating to this problem.  
Proposition \ref{PlD} gave a crude lower bound on the function $\ell_*(\Delta)$ defined as the infimum value of $\ell$ over all SIRSNs with the given value of $\Delta$. 
\begin{OP}
\label{OP2}
(i) Give quantitative estimates of the function $\ell_*(\Delta)$, improving Proposition \ref{PlD}. \\
(ii) Do ``optimal" networks attaining the infimum exist, and (if so) can we say something about the structure of the associated optimal networks?
\end{OP}
One  might make the (vague) conjecture that for some value of 
$\Delta$ the optimal network exploits 4-fold symmetry in some way analogous to our section \ref{sec-constr1} model, and that for some other value it
exploits 6-fold symmetry.

In section \ref{sec-Steiner} we showed (\ref{l*lb})  that the overall minimum normalized length $\ell_*:= \inf_\Delta   \ell_*(\Delta)$ 
satisfies $\ell_*\ge \sqrt{1/8}$.  
The third result was Lemma \ref{Lellp}, showng $\ell \le 2p(1)$.

\subsection{Traffic intensity}
\label{sec-traffic}
As mentioned in section \ref{sec-TS-2}, 
the conceptual point of $\EE(\infty,r)$ is to capture the idea of the 
major road - minor road spectrum, and the particular definition of $\EE(\infty,r)$
is mathematically convenient because of the scaling property (\ref{pr-scale})
of the edge-intensity $p(r)$.  
But from a real-world perspective it seems more natural to use 
some notion of traffic intensity.
Given any measure $\psi$ on source-destination pairs $(z_1,z_2)$, then
length measure $\Leb_1$
along
the routes $\RR(z_1,z_2)$ in a SIRSN induces a ``traffic intensity" measure 
$\widetilde{\psi}$ on $\cup_r \EE(\infty,r)$.
The natural measures $\psi$ to consider are specified by \\
(i) $z_1$ has Lebesgue measure $\Leb_2$ on $\Reals^2$\\
(ii) given $z_1$, the measure on $z:= z_2 - z_1$ 
has density $|z|^{-\beta}$.

\noindent 
The action of $\sigma_c$ on $\psi$ gives the measure specified by \\
(i) $z_1$ has measure $c^{-2} \Leb_2$ on $\Reals^2$\\
(ii) $z:= z_2 - z_1$  has density $c^{\beta -2} |z|^{- \beta}$\\
(iii) the measure along paths is $c^{-1} \Leb_1$\\
and this is the measure $\widetilde{\psi}$ scaled by $c^{\beta - 5}$.

To make a rigorous treatment, the issue is to show that $\widetilde{\psi}$ 
is a {\em locally finite} measure on $\EE(\infty,1)$.  
Heuristically one needs $\beta > 2$ so that the contribution from 
large $|z_2 - z_1|$ is finite, and $\beta < 4$ so that the contribution from 
small $|z_2 - z_1|$ is finite. 
\begin{OP}
Show that, perhaps under regularity assumptions on the SIRSN, 
for $2 < \beta < 4$ the construction above gives a locally finite measure 
 $\widetilde{\psi}$ on $\EE(\infty,1)$ 
 and hence on $\cup_r \EE(\infty,r)$.
\end{OP}

\subsection{Technical questions raised by results}

\subsubsection{Implications between different properties of a SIRSN} 
We have given various results of the form ``one property of a SIRSN implies another" 
for which we conjecture the converse is false.  
In particular, we expect there are counter-examples to most of the following, though of course this requires 
constructing  other examples of SIRSNs.

\begin{OP}
\label{OP-counter}
Prove, or give a counter-example to: \\
(i)   (\ref{ell-finite}) implies (\ref{p1intro}) \\
(ii) the unique singly-infinite geodesics property implies (\ref{prope01})\\
(iii) (\ref{prope01}) implies (\ref{def-strong-uniq}) \\
(iv) (\ref{def-strong-cont}) implies  (\ref{def-strong-uniq}).
\end{OP}

\subsubsection{Understanding the structure of $\EE(\infty, 1)$. }
We envisage $\EE(\infty, 1)$ as looking somewhat like a real-world network of major roads, but it is not clear
what aspects of real networks appear automatically in our SIRSN model. 
For instance, a priori $\EE(\infty, 1)$ need not be connected (it might contain a short segment in the middle of a route 
between two points at distance $2 + \eps$ apart) but it must 
contain an unbounded connected component (most of a singly-infinite geodesic).  
\begin{OP}
Does $\EE(\infty, 1)$ have a.s. only a single unbounded connected component?
\end{OP}

\subsubsection{Questions about lengths}
Even though we started the whole topic of SIRSNs by considering route-lengths, they have played a rather small role in our results, and many questions about route-lengths could be asked.
\begin{OP}
\label{OP-len}
Under what extra assumptions (if any) is it true that, for 
$U_1, U_2, \ldots$ independent uniform on $\disc(\origin, 1)$,
\[ \Ex \sup_{i \ge 1} \len [\RR(\origin, U_i)] < \infty ? \]
\end{OP}
The following (intuitively obvious) claim seems curiously hard to prove;   
 the difficulty lies in showing that the spanning subnetwork does not have (necessarily with low probability) 
huge length a long way away from the square.
\begin{OP}
Take $k$ uniform random points $Z_1,\ldots,Z_k$ in a square of area $k$ 
and consider the length $\len [  \span(Z_1,\ldots,Z_k) ]$ of the spanning subnetwork 
random network $\span(Z_1,\ldots,Z_k)$.
Prove
\[ \Ex \  \len [  \span(Z_1,\ldots,Z_k) ] \sim \ell k \mbox{ as } k \to \infty  . \]
\end{OP}

\subsection{Alternative starting points for a setup}

We started the whole modeling process by assuming we are given routes between points, but one can
imagine two different starting points.  
The first involves starting with a network of major roads and then adding successively more minor roads, so eventually the road network is dense in the plane.  
In other words, base a model on some {\em explicit} construction as $r$ decreases of some process $(\EE(r),\  \infty > r > 0)$ of ``roads of size $\ge r$" 
(in our setup this is achieved implicitly by the networks $\EE(\infty,r)$).  
Of course this corresponds to what we see when  zooming in on an online map of the real-world road network; 
the maps are designed to show only the relatively major roads 
within the window, and hence to show progressively more minor roads as one zooms in. 
In talks we show such zooms along with the online ``zooming in" demonstration \cite{wiki:BMscaling}  
of {\em Brownian scaling} to illustrate the concept of scale-invariance.  

The second, mathematically abstract, approach is to start with a random metric $d(z,z^\prime)$ on the plane, and define routes as geodesics.

But a technical difficulty with both of these approaches is that  there seems no simple way to guarantee 
{\em unique} routes between a.a. pairs of points in the plane -- in general one needs to add an {\em assumption} of uniqueness. 
The explicit models constructed in section \ref{sec-constr1} and outlined in 
section  \ref{sec:otherC} do use the ``random metric" idea, but the hard part of the construction is proving the uniqueness of routes, even in these simplest models we can imagine.
It is perhaps remarkable that our approach, taking routes as given with only the route-compatability property but with no explicit requirement 
that routes be minimum-cost in some sense, does lead to some non-obvious results.

\subsection{Empirical evidence of scale-invariance?}
For real-world road networks, can scale-invariance be even roughly true 
over some range of distance?
We mentioned one explicit piece of evidence (ordered segment lengths) 
in section \ref{sec-visualizing}; 
there is also evidence that mean route length is indeed 
roughly proportional to distance, 
though this is also consistent with other (non scale-invariant) models \cite{me-spatial-5}.

An interesting project would be to study the spanning subnetworks 
on (say) 4 real-world addresses, whose positions form roughly a square, 
randomly positioned, and find the empirical frequencies with 
which the various topologically different networks appear. 
Scale-invariance predicts these frequencies should not vary with the 
side-length of square; is this true?

\subsection{Other related literature}
\label{sec -ORL}

\subsubsection{Hop count in spatial networks}
\label{sec-hopcount}
There has been study of spatial networks with respect to the trade-off
between total network length
and average {\em graph distance} (hop count), instead of route-length.
See \cite{MR2510456} for a recent literature survey and empirical analysis.

\subsubsection{Continuum random trees in the plane}
\label{sec-tree}
Existence of continuum limits of discrete models of random trees has been conjectured, and studied non-rigorously in statistical physics,  for a long time, and since 2000 
spectacular  progress has been made on rigorous proofs.
For three models of random trees 
(uniform random spanning tree on $\Ints^2$,  minimal spanning tree on $\Ints^2$ (with random edge lengths), and the Euclidean minimal spanning tree on Poisson points),
\cite{MR1716768} established a rigorous ``tightness" result and gave sample properties of subsequential limits.  
A subsequent deep result \cite{MR2044671} established the existence of a continuum limit 
in the first model.
In these limits the paths have Hausdorff dimension greater than $1$ so
$D_1 = \infty$ a.s..
There should be a simple proof of the following, because our definition of SIRSN 
requires $\Ex D_1 = \infty$.
\begin{OP}
In a SIRSN, the subnetwork $\SS(1)$ cannot be a tree (with Steiner points).
\end{OP}

\subsubsection{Geodesics in first-passage percolation}
Geodesics in particular models of first-passage percolation have been studied 
in \cite{MR1387641}. It 
is unclear whether there is any substantial connection between the behavior of
geodesis in that setting and in our setting.

\subsubsection{A Monge-Kantorovitch approach}
\label{sec-MK}
A completely different approach to continuum networks, starting from
Monge-Kantorovitch optimal transport theory, is developed in the monograph 
by Buttazzo et al. \cite{MR2469110}.  Their model assumes\\
(i) some continuous distribution of sources and sinks \\
(ii) an {\em a priori} arbitrary set $\Sigma$ representing location of roads \\
(iii) two different costs-per-unit-length for travel inside [resp. outside] 
$\Sigma$.\\
An optimal network in one that minimizes total transportation cost for a given cost
functional on $\Sigma$.  
It is shown that, under regularity conditions, the 
optimal network is covered
by a finite number of Lipschitz curves of uniformly bounded length,
although it may have even uncountably many connected components.
But this theory does not seem to address statistics analogous to our $\Delta$ and 
$\ell$ in any quantitative way.

\subsubsection{The method of exchangeable substructures}

The general methodology of studying complicated random structures 
by studying induced substructures on random points has many applications \cite{me-kingman}.
In particular, the {\em Brownian continuum random tree} \cite{MR1166406} provides an
analogy 
for what we would like to see (Open Problem \ref{OP1})  in some
``mathematically natural" SIRSN -- see e.g. the formula (13) therein 
for the distribution of the induced subtree on random points -- though that is in the ``mean-field" setting 
without any $d$-dimensional geometry.

\subsubsection{Urban road networks.}
There is scattered literature on models for {\em urban} road networks, 
mostly with a rather different focus, though \cite{LGH06} has some conceptual similarities 
with our work.  

\subsubsection{Dynamic random graphs.} 
Conceptually, what we are doing with routes $\RR(z_1,z_2)$ and subnetworks $\SS(\lambda)$ is 
{\em exploring} a given network.  This is conceptually distinct from using sequential {\em constructions} of a network, a topic often called 
{\em dynamic random graphs} \cite{MR2271734}, even though the particular ``dynamic Gabriel" model outlined 
in section \ref{sec-OSM} does fit the ``dynamic" category.

\paragraph{Acknowldgements} 
My thanks to Justin Salez for details of the proof of Lemma \ref{Lccc}, to Wilfrid Kendall for ongoing collaboration, and to Cliff Stein for references to the algorithmic literature.

\appendix
\section{Appendix: A topology on the space of feasible subnetworks.}
We first define convergence of routes. 
Recall a feasible route $\rr(z,z^\prime)$ involves line segments between points 
$(z_i)$ which we will call {\em turn points} of the route. 
Given $\eps < |z^\prime - z|/2$ 
there is (starting from $z$) a last turn point 
$z_{(\eps)}$ before the route $\rr(z,z^\prime)$  first exits $\disc(z,\eps)$ 
and 
there is (starting from $z^\prime$) a last turn point 
$z^\prime_{(\eps)}$ before the reverse route $\rr(z^\prime,z)$  
first exits $\disc(z^\prime,\eps)$. 
Define
\[
\rr(z(n),z^\prime(n)) \to \rr(z,z^\prime)
\] 
 to mean \\
(i) $z(n) \to z, \ z^\prime (n) \to z^\prime \neq z$. \\
(ii) For each $\eps < |z^\prime - z|/2$ such that $\circl(z,\eps)$ and 
$\circl(z^\prime,\eps)$ do not contain any turn point of $\rr(z,z^\prime)$,
writing the turn points 
of the subroutes $\rr(z_{(\eps)}(n),z^\prime_{(\eps)}(n))$ 
as $(y_0(n), y_1(n),\ldots, y_k(n))$, we have 
\[ (y_0(n), y_1(n),\ldots, y_k(n)) \to (y_0, y_1,\ldots, y_k) \]
the limit being the turn points 
of the subroute $\rr(z_{(\eps)},z^\prime_{(\eps)})$, 
where $k$ is finite from the definition of feasible route. \\
(iii) The total lengths $L_{(\eps)}(n)$ of 
$\rr(z(n),z^\prime(n)) \cap(\disc(z,\eps) \cup \disc(z^\prime,\eps))$
satisfy 
\[ \lim_{\eps \to 0} \limsup_n L_{(\eps)}(n) = 0 . \]
Despite its inelegant formulation, this seems the ``natural" 
notion of convergence.

Now we specify, in a way analogous to (i-iii) above, 
what it means for a sequence $\ss(n)$ of feasible subnetworks 
on locally finite sets $\bz(n) = \{z^i(n)\}$ to converge to a limit subnetwork 
$\ss$ on $\bz$. \\ 
 (i) We need $\bz(n)$ to converge to $\bz$ in the usual sense of
convergence of simple point processes \cite{MR2371524}.  
This is equivalent to saying that if we take 
any $R$ such that $\circl(\origin,R)$ 
contains no point of $\bz$, then we can label the points of 
$\bz(n) \cap \disc(\origin,R)$ as $(z^1(n), \ldots, z^K(n))$ 
in such a way that 
$(z^1(n), \ldots, z^K(n)) \to (z^1, \ldots, z^K)$ , the limit 
(here and in analogous assertions below) being the points of $\bz \cap \disc(\origin,R)$. \\
(ii) Take $R$ and $(z^1, \ldots, z^K)$ as above and take
$\eps < \sfrac{1}{2} \min_{1\le i < j \le K} |z^i - z^j|$ 
such that $\cup_{1 \le i \le K} \circl(z^i,\eps)$ 
does not contain any turn point within $\ss$.  
Then we can label the turn points of 
$\cup_{1 \le i \le K} \rr(z^i_{(\eps)}(n), z^j_{(\eps)}(n))$ 
as $(y^u(n), 1 \le u \le L)$
in such a way that \\
 $(y^u(n), 1 \le u \le L)\to  (y^u, 1 \le u \le L)$ \\
 $(y^u(n),y^v(n))$ is an edge-segment of route $\rr(z^i(n),z^j(n))$ iff 
  $(y^u,y^v)$ is an edge-segment of route $\rr(z^i,z^j)$. \\
  (iii) For each $1 \le i < j \le K$ the routes  $\rr(z^i(n),z^j(n))$   satisfy (iii) above.
  \\
  (iv) Lemma \ref{Lss1} (i) implies that given $R$, the following quantity 
  (referring to the subnetwork $\ss$) is finite:
  \[ R^* := \min \{r: \
\bigcup_{z_i, z_j \in  \disc(\origin, r)}    \rr(z_i,z_j) \cap \disc(\origin, R)  =   \ss \cap \disc(\origin, R)  \   \} \] 
(that is, each edge of $\ss$ within $\disc(\origin, R)$ is part of some route between endpoints in $\disc(\origin, R^*)$).
We require 
\[ \limsup_n R^*(n) < \infty   \mbox{ for each } R < \infty . \]
 
We have described sequential convergence within the space of feasible subnetworks.  
It is routine to show this is convergence in some 
 complete separable metric space, but we won't pursue such theory here.

 \end{document}